\documentclass[12pt]{article}
\usepackage{amsmath}
\usepackage{amssymb}
\usepackage{vatola}

\def\q{\quad}
\def\qq{\qquad}
\def\qtq#1{\q\t{#1}\q}
\def\mod#1{\ (\text{\rm mod}\ #1)}
\def\t{\text}
\def\f{\frac}
\def\e{\equiv}
\def\b{\binom}
\def\sls#1#2{(\f{#1}{#2})}
 \def\ls#1#2{\big(\f{#1}{#2}\big)}
\def\Ls#1#2{\Big(\f{#1}{#2}\Big)}
\def\ap{\langle a\rangle_p}

\let \pro=\proclaim
\let \endpro=\endproclaim

\begin{document}
 \centerline {\bf
Congruences concerning binomial coefficients and binary quadratic
forms}
\par\q\newline
\centerline{Zhi-Hong Sun}\newline \centerline{School of Mathematics
and Statistics} \centerline{Huaiyin Normal University}
\centerline{Huaian, Jiangsu 223300, P.R. China} \centerline{Email:
zhsun@hytc.edu.cn} \centerline{Homepage:
http://maths.hytc.edu.cn/szh1.htm}
 \abstract{ \par Let $p>3$ be a prime. In this paper, we obtain the
 congruences for
 $$\sum_{k=0}^{p-1}\frac{w(k)\binom{2k}k^3}{(-8)^k},\
 \sum_{k=0}^{p-1}\frac{w(k)\binom{2k}k^2\binom{3k}k}{(-192)^k},\
 \sum_{k=0}^{p-1}\frac{w(k)\binom{2k}k^2\binom{4k}{2k}}{(-144)^k}\ \text{and}
 \ \sum_{k=0}^{p-1}\frac{w(k)\binom{2k}k^2\binom{4k}{2k}}{648^k}$$
modulo $p^2$, and partial results for $\sum_{k=0}^{(p-1)/2}
\binom{2k}k^3\frac{w(k)}{m^k}$ modulo $p^2$, where
$m\in\{1,16,-64,256,-512,4096\}$ and $w(k)\in\{k^2,k^3,\frac
1{k+1},\frac
 1{(k+1)^2},\frac 1{(k+1)^3},  \frac 1{2k-1},\frac 1{k+2}\}$.
 \par\q
 }
 \endabstract
\let\thefootnote\relax \footnotetext {The author was supported by
the National Natural Science Foundation of China (grant No.
12271200).}

\section*{1. Introduction}
\par\q The generalized binomial coefficient $\b ak$ is given by
$$\b a0=1\qtq{and}\b
ak=\f{a(a-1)\cdots(a-k+1)}{k!}\qtq{for}k=1,2,3,\ldots$$ It is easy
to see that $\b{-\f 12}k=\b{2k}k(-4)^{-k}$ and $p\mid \b {2k}k$ for
any odd prime $p$ satisfying $\f p2<k<p$. The binomial coefficients
$\b{2k}k$ $(k=1,2,3,\ldots)$ are called central binomial
coefficients. Set $C_k=\b{2k}k\f 1{k+1}$. Then $C_k\in\Bbb Z$ for
$k=0,1,2,\ldots$. The numbers $\{C_k\}$ are called Catalan numbers.
\par
 For positive integers $a,b$ and $n$, if
$n=ax^2+by^2$ for some integers $x$ and $y$, we briefly write that
$n=ax^2+by^2$. Let $p$ be an odd prime. It is known (see [BEW]) that
$$\align &p\e 1\mod 4\iff p=x^2+4y^2,
\\&p\e 1\mod 3\iff p=x^2+3y^2\q (p\not=3),
\\&p\e 1,2,4\mod 7\iff p=x^2+7y^2 \q(p\not=7),
\\&p\e 1,3\mod 8\iff p=x^2+2y^2.
\endalign$$
\par Let $p$ be a prime with $p\not=2,7$. In 1998, Ono[O] obtained the congruences
for $\sum_{k=0}^{p-1}\f 1{m^k}\b{2k}k^3$ modulo $p$ in the cases
$m=1,-8,16,-64,256,-512,4096$. For such values of $m$, in [Su1] the
author's brother Z.W. Sun conjectured the congruences for
$\sum_{k=0}^{p-1}\f 1{m^k}\b{2k}k^3$ modulo $p^2$, which have been
proved by the author in [S2] and Kibelbek et al in [KLMSY]. In [S7],
the author conjectured the congruences for $\sum_{k=0}^{p-1}\f
1{m^k}\b{2k}k^3$ modulo $p^3$. For instance, for any prime
$p\not=2,7$,
$$\sum_{k=0}^{p-1}\b{2k}k^3\e
\cases  4x^2-2p-\f{p^2}{4x^2}\mod {p^3}&\t{if $p\e 1,2,4\mod 7$ and
so $p=x^2+7y^2$,}
\\
-11p^2\b{[3p/7]}{[p/7]}^{-2} \mod {p^3}& \t{if $p\e 3\mod 7$,}
\\
-\f{11}{16}p^2\b{[3p/7]}{[p/7]}^{-2}\mod {p^3}& \t{if $p\e 5\mod
7$,}
\\
-\f{11}4p^2\b{[3p/7]}{[p/7]}^{-2}\mod {p^3}&\t{if $p\e 6\mod 7$,}
\endcases$$
where $[a]$ is the greatest integer not exceeding $a$. In [S7],
[S8], [S9] and [S10], the author posed many conjectures on the
congruences for $\sum_{k=0}^{(p-1)/2} \b{2k}k^3\f{w(k)}{m^k}$ modulo
$p^3$, where $m\in\{1,-8,16,64,$ $-64,256,-512,4096\}$ and
$$w(k)\in\Big\{k^2,k^3,\f 1{k+1},\f 1{(k+1)^2},\f 1{(k+1)^3}, \f
1{2k-1},\f 1{(2k-1)^2},\f 1{k+2},\f 1{k+3}\Big\}.$$ The case $m=64$
has been solved by the author in [S10]. We mention that Tauraso[T]
gave a congruence for $\sum_{k=0}^{p-1}\b{2k}k^3\f 1{64^k(k+1)^3}$
modulo $p^3$ via $p$-adic Gamma functions.
  As a typical conjecture in [S9],
for any prime $p\not=2,7$,
$$\sum_{k=0}^{(p-1)/2}\f{\b{2k}k^3}{k+1}
\e \cases -44y^2+2p\mod {p^3}&\t{if $p=x^2+7y^2\e 1,2,4\mod 7$,}
\\-\f{1}7\b{[3p/7]}{[p/7]}^2\mod p&\t{if $p\e 3\mod 7$,}
\\-\f{16}7\b{[3p/7]}{[p/7]}^2\mod p&\t{if $p\e 5\mod 7$,}
\\-\f{4}7\b{[3p/7]}{[p/7]}^2\mod p&\t{if $p\e 6\mod 7$.}
\endcases$$
\par It is known (see [S2-S5]) that
$$\b {-\f 13}k\b{-\f
23}k=\f{\b{2k}k\b{3k}k}{27^k}, \ \b {-\f 14}k\b{-\f
34}k=\f{\b{2k}k\b{4k}{2k}}{64^k}, \ \b {-\f 16}k\b{-\f
56}k=\f{\b{3k}{k}\b{6k}{3k}}{432^k}. \tag 1.1$$ For an odd prime $p$
let $\Bbb Z_p$ be the set of rational numbers whose denominator is
not divisible by $p$, and for an odd prime $p$ and $a\in\Bbb Z$ let
$\sls ap$ denote the Legendre symbol. Let $p>3$ be a prime. In 2003,
Rodriguez-Villegas[RV] made conjectures equivalent to
$$\align &\sum_{k=0}^{p-1}\f{\b{2k}k^2\b{3k}k}{108^k}\e
\cases 4x^2-2p\mod{p^2}&\t{if $p=x^2+3y^2\e 1\mod 3$,}
\\0\mod{p^2}&\t{if $p\e 2\mod 3$,}
\endcases
\\&\sum_{k=0}^{p-1}\f{\b{2k}k^2\b{4k}{2k}}{256^k}\e
\cases 4x^2-2p\mod{p^2}&\t{if $p=x^2+2y^2\e 1,3\mod 8$,}
\\0\mod{p^2}&\t{if $p\e 5,7\mod 8$,}
\endcases
\\&\Ls p3\sum_{k=0}^{p-1}\f{\b{2k}k\b{3k}k\b{6k}{3k}}{12^{3k}}\e
\cases 4x^2-2p\mod{p^2}&\t{if $p=x^2+4y^2\e 1\mod 4$,}
\\0\mod{p^2}&\t{if $p\e 3\mod 4$.}
\endcases\endalign$$
 These conjectures have been solved  by
Mortenson[M1] and Z.W. Sun[Su2]. In [Su1], [Su3] and [S1], Z.W. Sun
and the author posed many conjectures on the congruences for
$$\sum_{k=0}^{p-1}\f{\b{2k}k^2\b{3k}k}{m^k},\q
 \sum_{k=0}^{p-1}\f{\b{2k}k^2\b{4k}{2k}}{m^k}\qtq{and}
 \sum_{k=0}^{p-1}\f{\b{2k}k\b{3k}k\b{6k}{3k}}{m^k}$$
 modulo $p^2$, where $m\in\Bbb Z_p$ and $m\not\e 0\mod p.$
 In [S7], [S8], [S9] and [S10], the author posed further conjectures on the congruences for
$$\sum_{k=0}^{p-1}\f{w(k)\b{2k}k^2\b{3k}k}{m^k},\q
 \sum_{k=0}^{p-1}\f{w(k)\b{2k}k^2\b{4k}{2k}}{m^k}\qtq{and}
 \sum_{k=0}^{p-1}\f{w(k)\b{2k}k\b{3k}k\b{6k}{3k}}{m^k}$$
 modulo $p^3$, where $w(k)\in\{1,k,k^2,k^3,\f 1{k+1},\f
 1{(k+1)^2},\f 1{(k+1)^3},
 \f 1{2k-1},\f 1{k+2}\}$.

\par  Let $p>3$ be a prime and $a,m\in\Bbb Z_p$ with $m\not\e 0\mod p$. In [S5],
the author investigated the congruences concerning $\sum_{k=0}^{p-1}
\b ak\b{-1-a}k\f 1{m^k}$ modulo $p^2$. In this paper, we deduce the
congruences for $\sum_{k=0}^{p-1} w(k)\b{2k}k\b ak\b{-1-a}k\f
1{m^k}$ modulo $p^2$, where $w(k)\in\{k,k^2,k^3,\f 1{k+1},\f
 1{(k+1)^2},\f 1{(k+1)^3},
 \f 1{2k-1},\f 1{k+2}\}$. In particular,
 we  obtain the congruences for
$$\sum_{k=0}^{p-1}\f{w(k)\b{2k}k^3}{(-8)^k},\q
\sum_{k=0}^{p-1}\f{w(k)\b{2k}k^2\b{3k}k}{(-192)^k},\q
\sum_{k=0}^{p-1}\f{w(k)\b{2k}k^2\b{4k}{2k}}{(-144)^k}\qtq{and}
\sum_{k=0}^{p-1}\f{w(k)\b{2k}k^2\b{4k}{2k}}{648^k}$$ modulo $p^2$,
and partial results for $\sum_{k=0}^{(p-1)/2}
\b{2k}k^3\f{w(k)}{m^k}$ modulo $p^2,$ where $m\in\{1,16,-64,256,$
$-512,4096\}$. As a typical example, for any odd prime $p\e
1,2,4\mod 7$ and so $p=x^2+7y^2$ we show that
$$\align&\sum_{k=0}^{(p-1)/2}\f{\b{2k}k^3}{k+1}
\e -44y^2+2p\mod {p^2},
\\&\sum_{k=0}^{(p-1)/2}\f{\b{2k}k^3}{(k+1)^3}\e \f 18-\f{201}2y^2-\f
94p\mod {p^2},
\\&\sum_{k=0}^{(p-1)/2}k^3\b{2k}k^3\e
-\f{5408x^2}{27783}+\f{2992p}{9261}\mod {p^2}.
\endalign$$
 Also, for any prime $p>3$,
$$\sum_{k=0}^{(p-1)/2}\f{\b{2k}k^3}{(-8)^k(k+2)}
 \e \cases\f{64}{27}x^2-\f{35}{27}p\mod{p^2}&\t{if $p=x^2+4y^2\e 1\mod
 4$,}\\\f p9\mod {p^2}&\t{if $p\e 3\mod 4$.}
 \endcases$$
\par Throughout this paper, for an odd prime $p$ let
$$\align &R_1(p)=(2p+2-2^{p-1})\b{(p-1)/2}{[p/4]}^2,
\\&R_3(p)=\Big(1+2p+\f
43(2^{p-1}-1)-\f 32(3^{p-1}-1)\Big) \b{(p-1)/2}{[p/6]}^2.
\endalign$$

\section*{2. Congruences for $\sum_{k=0}^{p-2}
\b ak\b{-1-a}k\b{2k}k\f 1{m^k(k+1)}$
 modulo $p^2$}
\pro{Lemma 2.1} For $n=0,1,2,\ldots$ we have
$$\align&\sum_{k=0}^nk\b ak\b{-1-a}k(n+1-k)\b a{n+1-k}\b{-1-a}{n+1-k}
\\&=a(a+1)\sum_{k=0}^n\b
ak\b{-1-a}k\b{2k}{k+1}(-1)^{n+1-k}\b{k-1}{n-k}.\endalign$$
\endpro
Proof. For $n=0,1,2,\ldots$ set
$$\align&S(n)=\sum_{k=0}^nk\b ak\b{-1-a}k(n+1-k)\b
a{n+1-k}\b{-1-a}{n+1-k},
\\&T(n)=a(a+1)\sum_{k=0}^n\b
ak\b{-1-a}k\b{2k}{k+1}(-1)^{n+1-k}\b{k-1}{n-k}.
\endalign$$
Using the package sumtools in Maple, we find that for $n\ge 2$,
$$\align&(n^3-n)S(n)=2(n^3-(2a^2+2a+1)n+a(a+1))S(n-1)-(n^3-(2a+1)^2n)S(n-2),
\\&(n^3-n)T(n)=2(n^3-(2a^2+2a+1)n+a(a+1))T(n-1)-(n^3-(2a+1)^2n)T(n-2).
\endalign$$
This shows that $\{S(n)\}$ and $\{T(n)\}$ satisfy the same
recurrence relation. Also, $S(0)=T(0)=0$ and $S(1)=T(1)=a^2(a+1)^2$.
Thus, we must have $S(n)=T(n)$ for $n=0,1,2,\ldots$. This proves the
lemma. \vskip 0.25cm
\par {\bf Remark 2.1} By Lemma 2.1, we have
 $$\align&\Big(\sum_{k=0}^{\infty}\b ak\b{-1-a}kk(-t)^k\Big)^2
\\&=\sum_{n=0}^{\infty}\Big(\sum_{k=0}^n k\b ak\b{-1-a}k(n+1-k)\b
a{n+1-k}\b{-1-a}{n+1-k}\Big)(-t)^{n+1}
\\&=a(a+1)\sum_{n=0}^{\infty}\Big(\sum_{k=0}^n\b{2k}{k+1}\b
ak\b{-1-a}k(-1)^k\b{k-1}{n-k}\Big)t^{n+1}
\\&=a(a+1)\sum_{k=1}^{\infty}\b{2k}{k+1}\b
ak\b{-1-a}k(-1)^kt^{k+1}\sum_{r=0}^{k-1}\b{k-1}rt^r
\\&=\f{a(a+1)t}{t+1}\sum_{k=0}^{\infty}\b{2k}{k+1}\b
ak\b{-1-a}k(-t(t+1))^k.
\endalign$$
\vskip 0.2cm

\pro{Lemma 2.2 ([S6, Theorem 2.2])} Let $p$ be an odd prime and
$a\in\Bbb Z_p$. Then
$$\Big(\sum_{k=0}^{p-1}\b ak\b {-1-a}k(-t)^k\Big)^2\e
\sum_{k=0}^{p-1}\b{2k}k\b ak\b {-1-a}k(-t(t+1))^k\mod{p^2}.$$
\endpro
\par For the special cases $a=-\f 12,-\f 13,-\f 14,-\f 16$ see [S2], [S3]
and [S4].
 \pro{Theorem 2.1} Let $p$ be an odd prime, $a\in\Bbb Z_p$
and $a\not\e 0,-1\mod p$. Then
$$\align&\sum_{k=0}^{p-2}\b ak\b{-1-a}k\b{2k}k\f{(-t(t+1))^k}{k+1}
\\&\e\Big(\sum_{k=0}^{p-1}\b ak\b {-1-a}k(-t)^k\Big)^2
-\f{t+1}{a(a+1)t}\Big(\sum_{k=0}^{p-1}\b ak\b{-1-a}kk(-t)^k\Big)^2
\mod {p^2}.\endalign$$
\endpro
Proof. For $k\in\{\f{p+1}2,\ldots,p-2\}$ we see that $p\mid
\b{2k}{k+1}$ and
$$\b ak\b {-1-a}k=(-1)^k\f{(a+k)(a+k-1)\cdots (a-k+1)}{k!^2}\e 0\mod p.$$
Since $a\not\e 0,-1\mod p$, we also have $\b a{p-1}\e
\b{-1-a}{p-1}\e 0\mod p$. Thus, applying Lemma 2.1 we deduce that
$$\align&\f 1{t+1}\sum_{k=0}^{p-1}\b{2k}{k+1}\b ak\b {-1-a}k(-t(t+1))^k
\\&\e \sum_{k=1}^{(p-1)/2}\b{2k}{k+1}\b ak\b
{-1-a}k(-t)^k(t+1)^{k-1}
\\&=\sum_{k=1}^{(p-1)/2}\b{2k}{k+1}\b ak
\b {-1-a}k(-t)^k\sum_{r=0}^{k-1}\b {k-1}rt^r
\\&=\sum_{n=1}^{p-2}t^n
\sum_{k=1}^{\min\{n,\f{p-1}2\}}\b{2k}{k+1}\b ak\b{-1-a}k\b
{k-1}{n-k}(-1)^{k}
\\&\e\sum_{n=1}^{p-2}(-t)^n
\sum_{k=0}^n\b{2k}{k+1}\b ak\b{-1-a}k\b {k-1}{n-k}(-1)^{n-k}
\\&=-\f 1{a(a+1)}\sum_{n=1}^{p-2}(-t)^n
\sum_{k=0}^n\b ak\b{-1-a}kk(n+1-k)\b a{n+1-k}\b{-1-a}{n+1-k}
\\&=\f 1{a(a+1)t}\sum_{k=0}^{p-2}\b ak\b{-1-a}kk(-t)^k\sum_{n=k}^{p-2}
\b a{n+1-k}\b{-1-a}{n+1-k}(n+1-k)(-t)^{n+1-k}
\\&=\f 1{a(a+1)t}\sum_{k=0}^{p-2}\b ak\b{-1-a}kk(-t)^k
\\&\q\times\Big(\sum_{r=0}^{p-2} \b ar\b{-1-a}rr(-t)^r-\sum_{r=p-k}^{p-2} \b
ar\b{-1-a}rr(-t)^r\Big)\mod {p^2}.
\endalign$$
 Let $\ap$ be the least nonnegative residue of $a$ modulo $p$. For $\ap<k\le p-1$ we see that $\b ak\e 0\mod p$, and for $0\le
k\le \ap$ and $p-k\le r\le p-1$ we see that $r\ge p-k>p-1-\ap$ and
so $ \b{-1-a}r\e 0\mod p$. Hence, for $0\le k\le p-1$ and $p-k\le
r\le p-1$ we have $\b ak\b{-1-a}r\e 0\mod p$. Replacing $a$ with
$-1-a$ we get $\b{-1-a}k\b ar\e 0\mod p$ and so $ \b ak\b{-1-a}k\b
ar\b{-1-a}r\e 0\mod {p^2}$.  Thus,
$$\aligned&\f 1{t+1}\sum_{k=0}^{p-1}\b{2k}{k+1}\b ak\b {-1-a}k(-t(t+1))^k
\\&\e\f 1{a(a+1)t}\Big(\sum_{k=0}^{p-2}\b ak\b{-1-a}kk(-t)^k\Big)^2\mod
{p^2}.\endaligned\tag 2.1$$ Since $\b{2k}k\f
1{k+1}=\b{2k}k-\b{2k}{k+1}$ and $\b a{p-1}\e\b{-1-a}{p-1}\e
\b{2(p-1)}{p-1}\e 0\mod p$, we then get
$$\align&\sum_{k=0}^{p-2}\b ak\b{-1-a}k\b{2k}k\f{(-t(t+1))^k}{k+1}
\\&\e \sum_{k=0}^{p-1}\b ak\b
{-1-a}k\b{2k}{k}(-t(t+1))^k-\f{t+1}{a(a+1)t}\Big(\sum_{k=0}^{p-1}\b
ak\b{-1-a}kk(-t)^k\Big)^2\mod {p^2}.\endalign$$ Now applying Lemma
2.2 yields the result.
 \vskip 0.25cm
 \par Taking $a=-\f 12,-\f 13,-\f 14,-\f 16$ in Theorem 2.1 and then applying (1.1)
 yields the following corollaries.
 \pro{Corollary 2.1} Let $p$ be an odd prime. Then
$$\align&\sum_{k=0}^{(p-1)/2}\f{\b{2k}k^3}{k+1}\Big(-\f{t(t+1)}{16}\Big)^k
\\&\e \Big(\sum_{k=0}^{(p-1)/2}\b{2k}k^2\big(-\f t{16}\big)^k\Big)^2
+\f {4(t+1)}t\Big(\sum_{k=0}^{(p-1)/2}\b{2k}k^2k\big(-\f
t{16}\big)^k\Big)^2\mod {p^2}.\endalign$$
\endpro
\pro{Corollary 2.2} Let $p>3$ be a prime. Then
$$\align&\sum_{k=0}^{p-2}\f{\b{2k}k^2\b{3k}k}{k+1}\Big(-\f{t(t+1)}{27}
\Big)^k
\\&\e \Big(\sum_{k=0}^{p-1}\b{2k}k\b{3k}k\big(-\f t{27}\big)^k\Big)^2
+\f {9(t+1)}{2t}\Big(\sum_{k=0}^{p-1}\b{2k}k\b{3k}kk\big(-\f
t{27}\big)^k\Big)^2\mod {p^2}.\endalign$$
\endpro
\pro{Corollary 2.3} Let $p>3$ be a prime. Then
$$\align&\sum_{k=0}^{p-2}\f{\b{2k}k^2\b{4k}{2k}}{k+1}\Big(-\f{t(t+1)}{64}
\Big)^k
\\&\e \Big(\sum_{k=0}^{p-1}\b{2k}k\b{4k}{2k}\big(-\f t{64}\big)^k\Big)^2
+\f {16(t+1)}{3t}\Big(\sum_{k=0}^{p-1}\b{2k}k\b{4k}{2k}k\big(-\f
t{64}\big)^k\Big)^2\mod {p^2}.\endalign$$
\endpro
\pro{Corollary 2.4} Let $p>5$ be a prime. Then
$$\align&\sum_{k=0}^{p-2}\f{\b{2k}k\b{3k}k\b{6k}{3k}}{k+1}
\Big(-\f{t(t+1)}{432} \Big)^k
\\&\e \Big(\sum_{k=0}^{p-1}\b{3k}k\b{6k}{3k}\big(-\f t{432}\big)^k\Big)^2
+\f {36(t+1)}{5t}\Big(\sum_{k=0}^{p-1}\b{3k}k\b{6k}{3k}k\big(-\f
t{432}\big)^k\Big)^2\mod {p^2}.\endalign$$
\endpro

\pro{Theorem 2.2} Let $p$ be an odd prime, $a\in\Bbb Z_p$ and
$a\not\e 0,-1\mod p$, and let $t$ be a $p$-adic integer such that
$t(t+1)\not\e 0\mod p$. If $\sum_{k=0}^{p-1}\b
ak\b{-1-a}k\b{2k}k(-t(t+1))^k\not\e 0\mod p$, then
$$\align \sum_{k=0}^{p-2}&\f{\b ak\b{-1-a}k\b{2k}k}{k+1}(-t(t+1))^k
\e \sum_{k=0}^{p-1}\b ak\b{-1-a}k\b{2k}k(-t(t+1))^k
\\&\q\qq\q-\f{(2t+1)^2}{4a(a+1)t(t+1)}\cdot \f{(\sum_{k=0}^{p-1}\b
ak\b{-1-a}k\b{2k}kk(-t(t+1))^k)^2}{\sum_{k=0}^{p-1}\b
ak\b{-1-a}k\b{2k}k(-t(t+1))^k}\mod {p^2}.\endalign$$
\endpro
Proof. By Lemma 2.2,
$$\Big(\sum_{k=0}^{p-1}\b ak\b {-1-a}k(-t)^k\Big)^2=
\sum_{k=0}^{p-1}\b{2k}k\b ak\b {-1-a}k(-t(t+1))^k+p^2f(t),$$ where
$f(t)$ is a polynomial in $t$ with coefficients in $\Bbb Z_p$.
Taking derivatives on both sides yields
$$\align&2\Big(\sum_{k=0}^{p-1}\b ak\b {-1-a}k(-t)^k\Big)
\Big(-\sum_{k=0}^{p-1}\b ak\b {-1-a}kk(-t)^{k-1}\Big) \\&=
\sum_{k=0}^{p-1}\b{2k}k\b ak\b
{-1-a}kk(-t(t+1))^{k-1}(-2t-1)+p^2f'(t).\endalign$$ This implies
that
$$\aligned &2\Big(\sum_{k=0}^{p-1}\b ak\b {-1-a}k(-t)^k\Big)
\Big(\sum_{k=0}^{p-1}\b ak\b {-1-a}kk(-t)^k\Big)
\\&\e \f{2t+1}{t+1}
\Big(\sum_{k=0}^{p-1}\b{2k}k\b ak\b {-1-a}kk(-t(t+1))^k\Big) \mod
{p^2}.\endaligned\tag 2.2$$ From (2.2) and Lemma 2.2 we deduce that
$$\aligned &
\Big(\sum_{k=0}^{p-1}\b ak\b
{-1-a}kk(-t)^k\Big)^2\Big(\sum_{k=0}^{p-1}\b{2k}k\b ak\b
{-1-a}k(-t(t+1))^k\Big)
\\&\e \f{(2t+1)^2}{4(t+1)^2}
\Big(\sum_{k=0}^{p-1}\b{2k}k\b ak\b {-1-a}kk(-t(t+1))^k\Big)^2 \mod
{p^2}.\endaligned\tag 2.3$$ If $\sum_{k=0}^{p-1}\b{2k}k\b ak\b
{-1-a}k(-t(t+1))^k\not\e 0\mod p$, appealing to Lemma 2.2, Theorem
2.1 and (2.3) we derive that
$$\align& \sum_{k=0}^{p-2}\f{\b ak\b{-1-a}k\b{2k}k}{k+1}(-t(t+1))^k
-\sum_{k=0}^{p-1}\b ak\b{-1-a}k\b{2k}k(-t(t+1))^k
\\&\e -\f{t+1}{a(a+1)t}\Big(\sum_{k=0}^{p-1}\b ak\b
{-1-a}kk(-t)^k\Big)^2
\\&\e -\f{t+1}{a(a+1)t}\cdot\f{(2t+1)^2}{4(t+1)^2}
\cdot \f{(\sum_{k=0}^{p-1}\b
ak\b{-1-a}k\b{2k}kk(-t(t+1))^k)^2}{\sum_{k=0}^{p-1}\b
ak\b{-1-a}k\b{2k}k(-t(t+1))^k}\mod {p^2}.\endalign$$ This yields the
result.

\vskip 0.25cm \par Taking $a=-\f 12,-\f 13,-\f 14,-\f 16$ in Theorem
2.2 and then applying (1.1) gives the following corollaries.

\pro{Corollary 2.5} Let $p$ be an odd prime. If $t$ is a $p$-adic
integer such that
$t(t+1)\sum_{k=0}^{(p-1)/2}\b{2k}k^3(-\f{t(t+1)}{16})^k\not\e 0\mod
p$, then
$$\align &\sum_{k=0}^{(p-1)/2}\f{\b{2k}k^3}{k+1}\Big(-\f{t(t+1)}{16}\Big)^k
\\&\e \sum_{k=0}^{(p-1)/2}\b{2k}k^3\Big(-\f{t(t+1)}{16}\Big)^k
+\f{(2t+1)^2}{t(t+1)}\cdot\f{(\sum_{k=0}^{(p-1)/2}\b{2k}k^3k(-\f{t(t+1)}{16})^k)^2}
{\sum_{k=0}^{(p-1)/2}\b{2k}k^3(-\f{t(t+1)}{16})^k}\mod {p^2}.
\endalign$$
\endpro

\pro{Corollary 2.6} Let $p>3$ be a prime. If $t$ is a $p$-adic
integer such that
$t(t+1)\sum_{k=0}^{p-1}\b{2k}k^2\b{3k}{k}(-\f{t(t+1)}{27})^k\not\e
0\mod p$, then
$$\align&\sum_{k=0}^{p-2}\f{\b{2k}k^2\b{3k}{k}}{k+1}\Big(-\f{t(t+1)}{27}
\Big)^k
\\&\e \sum_{k=0}^{p-1}\b{2k}k^2\b{3k}{k}\Big(-\f{t(t+1)}{27}\Big)^k
+\f {9(2t+1)^2}{8t(t+1)}\cdot
\f{(\sum_{k=0}^{p-1}\b{2k}k^2\b{3k}{k}k(-\f{t(t+1)}{27})^k)^2}
{\sum_{k=0}^{p-1}\b{2k}k^2\b{3k}{k}(-\f{t(t+1)}{27})^k}\mod {p^2}.
\endalign$$
\endpro
\pro{Corollary 2.7} Let $p>3$ be a prime. If $t$ is a $p$-adic
integer such that
$t(t+1)\sum_{k=0}^{p-1}\b{2k}k^2\b{4k}{2k}(-\f{t(t+1)}{64})^k\not\e
0\mod p$, then
$$\align&\sum_{k=0}^{p-2}\f{\b{2k}k^2\b{4k}{2k}}{k+1}\Big(-\f{t(t+1)}{64}\Big)^k
\\&\e \sum_{k=0}^{p-1}\b{2k}k^2\b{4k}{2k}\Big(-\f{t(t+1)}{64}\Big)^k
+\f {4(2t+1)^2}{3t(t+1)}\cdot
\f{(\sum_{k=0}^{p-1}\b{2k}k^2\b{4k}{2k}k(-\f{t(t+1)}{64})^k)^2}
{\sum_{k=0}^{p-1}\b{2k}k^2\b{4k}{2k}(-\f{t(t+1)}{64})^k}\mod {p^2}.
\endalign$$
\endpro
 \pro{Corollary 2.8} Let $p>5$ be a prime. If $t$ is a $p$-adic
integer such that
$t(t+1)\sum_{k=0}^{p-1}\b{2k}k\b{3k}{k}\b{6k}{3k}(-\f{t(t+1)}{432})^k\not\e
0\mod p$, then
$$\align\sum_{k=0}^{p-2}&\f{\b{2k}k\b{3k}{k}\b{6k}{3k}}{k+1}
\Big(-\f{t(t+1)}{432} \Big)^k \e
\sum_{k=0}^{p-1}\b{2k}k\b{3k}{k}\b{6k}{3k}\Big(-\f{t(t+1)}{432}\Big)^k
\\&\q+\f {9(2t+1)^2}{5t(t+1)}\cdot
\f{(\sum_{k=0}^{p-1}\b{2k}k\b{3k}{k}\b{6k}{3k}k(-\f{t(t+1)}{432})^k)^2}
{\sum_{k=0}^{p-1}\b{2k}k\b{3k}{k}\b{6k}{3k}(-\f{t(t+1)}{432})^k}\mod
{p^2}.
\endalign$$
\endpro

 \pro{Theorem 2.3} Let $p$ be an odd prime such that $p\e 1,2,4\mod
7$ and so $p=x^2+7y^2$. Then
$$ \align &\sum_{k=0}^{(p-1)/2}\f{\b{2k}k^3}{k+1}\e -44y^2+2p\mod
{p^2},
\\& \sum_{k=0}^{(p-1)/2}\f{\b{2k}k^3}{4096^k(k+1)}\e
(-1)^{\f{p-1}2}(72y^2+2p)\mod {p^2}.\endalign$$
\endpro
Proof. In [Su1], Z.W. Sun showed that
$$\sum_{k=0}^{p-1}(21k+8)\b{2k}k^3\e 8p\mod {p^3}.$$
By [S2, Theorem 3.4],
$$\sum_{k=0}^{(p-1)/2}\b{2k}k^3\e 4x^2-2p=-28y^2+2p \mod {p^2}.
\tag 2.4$$ Thus, $$\sum_{k=0}^{(p-1)/2}k\b{2k}k^3\e \f
8{21}(3p-4x^2)=\f 8{21}(28y^2-p)\mod {p^2}.\tag 2.5$$ Set
$t=\f{3\sqrt{-7}-1}2$. Then $t(t+1)=-16$ and $(2t+1)^2=-63$. Now,
from Corollary 2.5 and the above we deduce that
$$\align \sum_{k=0}^{(p-1)/2}\f{\b{2k}k^3}{k+1}
&\e\sum_{k=0}^{(p-1)/2}\b{2k}k^3+\f{63}{16}\cdot
\f{(\sum_{k=0}^{(p-1)/2}k\b{2k}k^3)^2}{\sum_{k=0}^{(p-1)/2}\b{2k}k^3}
\\&\e
-28y^2+2p+\f{63}{16}\cdot\f{8^2}{21^2}\cdot\f{(28y^2-p)^2}{-28y^2+2p}
\\&\e -28y^2+2p-16\cdot\f{28y^4-2py^2}{28y^2-2p}=-44y^2+2p
\mod {p^2}.\endalign$$
\par Let us consider the remaining part.  From
[S2, Theorem 3.4],
$$\sum_{k=0}^{(p-1)/2}\f{\b{2k}k^3}{4096^k}\e (-1)^{\f{p-1}2}(4x^2-2p)
=(-1)^{\f{p-1}2}(-28y^2+2p) \mod {p^2}.\tag 2.6$$ By [OZ],
$$\sum_{k=0}^{(p-1)/2}(42k+5)\f{\b{2k}k^3}{4096^k}\e
5(-1)^{\f{p-1}2}p\mod {p^2}.$$ Thus,
$$ \sum_{k=0}^{(p-1)/2}\f{k\b{2k}k^3}{4096^k}
\e \f
1{42}\Big(5(-1)^{\f{p-1}2}p-5\sum_{k=0}^{(p-1)/2}\f{\b{2k}k^3}{4096^k}\Big)
\e \f 5{42}(-1)^{\f{p-1}2}(28y^2-p)\mod {p^2}.\tag 2.7 $$ Set $t=\f
12(\f {3\sqrt {7}}8-1)$. Then $t(t+1)=-\f 1{256}$ and
$(2t+1)^2=\f{63}{64}$. From the above and Corollary 2.5 we deduce
that
$$\align \sum_{k=0}^{(p-1)/2}\f{\b{2k}k^3}{4096^k(k+1)}
&\e \sum_{k=0}^{(p-1)/2}\f{\b{2k}k^3}{4096^k}
-\f{63/64}{1/256}\cdot\f{(\sum_{k=0}^{(p-1)/2}k\b{2k}k^3
4096^{-k})^2} {\sum_{k=0}^{(p-1)/2}\b{2k}k^34096^{-k}}
\\&\e (-1)^{\f{p-1}2}(-28y^2+2p)-63\cdot 4\cdot\f {5^2}{42^2}
\cdot \f{(28y^2-p)^2}{(-1)^{\f{p-1}2}(-28y^2+2p)}
\\&\e (-1)^{\f{p-1}2}(72y^2+2p)\mod {p^2}.
\endalign$$
This completes the proof. \vskip 0.25cm
 \pro{Theorem 2.4} Let $p$ be a prime such
that $p\e 1\mod 3$ and so $p=x^2+3y^2$. Then
$$ \align &\sum_{k=0}^{(p-1)/2}\f{\b{2k}k^3}{16^k(k+1)}\e -16y^2+2p\mod
{p^2},
\\& \sum_{k=0}^{(p-1)/2}\f{\b{2k}k^3}{256^k(k+1)}\e
(-1)^{\f{p-1}2}(-8y^2+2p)\mod {p^2}.\endalign$$
\endpro
Proof. By [S2, Theorem 3.4],
$$\sum_{k=0}^{(p-1)/2}\f{\b{2k}k^3}{16^k}\e 4x^2-2p\e
-12y^2+2p\mod {p^2}.\tag 2.8$$ By [He],
$$\sum_{k=0}^{(p-1)/2}\f{k\b{2k}k^3}{16^k}\e p-\f 43x^2=-\f
p3+4y^2\mod {p^2}.\tag 2.9$$ Set $t=(\sqrt{-3}-1)/2$. Then
$t(t+1)=-1$ and $(2t+1)^2=-3$. Now, from the above and Corollary 2.5
we deduce that $$\align \sum_{k=0}^{(p-1)/2}\f{\b{2k}k^3}{16^k(k+1)}
&\e \sum_{k=0}^{(p-1)/2}\f{\b{2k}k^3}{16^k}
+3\cdot\f{(\sum_{k=0}^{(p-1)/2}k\b{2k}k^3 16^{-k})^2}
{\sum_{k=0}^{(p-1)/2}\b{2k}k^316^{-k}}
\\&\e -12y^2+2p+3\cdot\f{(4y^2-\f p3)^2}{-12y^2+2p}\e
-12y^2+2p+\f{48y^4-8py^2}{-12y^2+2p}\\&=-16y^2+2p\mod {p^2}.
\endalign$$
This proves the first part.
\par Now we consider the remaining part. By [S2, Theorem 3.4],
$$\sum_{k=0}^{(p-1)/2}\f{\b{2k}k^3}{256^k}\e
(-1)^{\f{p-1}2}(4x^2-2p) =(-1)^{\f{p-1}2}(-12y^2+2p)\mod {p^2}.\tag
2.10$$ By [L],
$$\sum_{k=0}^{(p-1)/2}(6k+1)\f{\b{2k}k^3}{256^k}\e (-1)^{\f{p-1}2}p
\mod {p^4}.$$ Thus,
$$\sum_{k=0}^{(p-1)/2}k\f{\b{2k}k^3}{256^k}
\e \f 16\Big((-1)^{\f{p-1}2}p-
\sum_{k=0}^{(p-1)/2}\f{\b{2k}k^3}{256^k}\Big)\e \f
16(-1)^{\f{p-1}2}(12y^2-p)\mod {p^2}.\tag 2.11$$ Set $t=\f{\sqrt
3-2}4$. Then $t(t+1)=-\f 1{16}$ and $(2t+1)^2=\f 34$. Now, from the
above and Corollary 2.5 we derive that
$$\align\sum_{k=0}^{(p-1)/2}\f{\b{2k}k^3}{256^k(k+1)}
&\e \sum_{k=0}^{(p-1)/2}k\f{\b{2k}k^3}{256^k}-\f{3/4}{1/16}\cdot
\f{(\sum_{k=0}^{(p-1)/2}k\b{2k}k^3 256^{-k})^2}
{\sum_{k=0}^{(p-1)/2}\b{2k}k^3256^{-k}}
\\&\e
(-1)^{\f{p-1}2}(-12y^2+2p)-\f{12}{36}\cdot\f{(12y^2-p)^2}
{(-1)^{\f{p-1}2}(-12y^2+2p)} \\&\e (-1)^{\f{p-1}2}(-8y^2+2p) \mod
{p^2}.\endalign$$ This completes the proof. \vskip 0.25cm
\par\q
\pro{Theorem 2.5} Let $p$ be a prime such that $p\e 1,3\mod 8$ and
so $p=x^2+2y^2$. Then
$$\sum_{k=0}^{(p-1)/2}\f{\b{2k}k^3}{(-64)^k(k+1)}
\e (-1)^{\f{p-1}2}(-12y^2+2p)\mod{p^2}.$$
\endpro
Proof. By [S2, Theorem 3.4],
$$\sum_{k=0}^{(p-1)/2}\f{\b{2k}k^3}{(-64)^k}\e
(-1)^{\f{p-1}2}(4x^2-2p)=(-1)^{\f{p-1}2}(-8y^2+2p)\mod {p^2}.\tag
2.12$$ By [M2],
$$\sum_{k=0}^{(p-1)/2}(4k+1)\f{\b{2k}k^3}{(-64)^k}\e
(-1)^{\f{p-1}2}p\mod {p^3}.$$ Thus,
$$\sum_{k=0}^{(p-1)/2}\f{k\b{2k}k^3}{(-64)^k}
\e \f
14\Big((-1)^{\f{p-1}2}p-\sum_{k=0}^{(p-1)/2}\f{\b{2k}k^3}{(-64)^k}\Big)
\e \f 14(-1)^{\f{p-1}2}(8y^2-p)\mod {p^2}.\tag 2.13$$ Set
$t=\f{\sqrt 2-1}2$. Then $t(t+1)=\f 14$ and $(2t+1)^2=2$. Now, from
the above and Corollary 2.5 we deduce that
$$\align\sum_{k=0}^{(p-1)/2}\f{\b{2k}k^3}{(-64)^k(k+1)}
&\e \sum_{k=0}^{(p-1)/2}k\f{\b{2k}k^3}{(-64)^k}+\f{2}{1/4}\cdot
\f{(\sum_{k=0}^{(p-1)/2}k\b{2k}k^3(-64)^{-k})^2}
{\sum_{k=0}^{(p-1)/2}\b{2k}k^3(-64)^{-k}}
\\&\e (-1)^{\f{p-1}2}(-8y^2+2p)+\f 8{16}\cdot\f{(8y^2-p)^2}{(-1)^{\f{p-1}4}(-8y^2+2p)}
\\&\e (-1)^{\f{p-1}2}(-12y^2+2p)\mod {p^2}.\endalign$$
This proves the theorem. \vskip 0.25cm
 \pro{Theorem 2.6} Let $p$ be
a prime such that $p\e 1\mod 4$ and so $p=x^2+4y^2$. Then
$$\sum_{k=0}^{(p-1)/2}\f{\b{2k}k^3}{(-512)^k(k+1)} \e
(-1)^{\f{p-1}4}(-32y^2+2p)\mod{p^2}.$$
\endpro
Proof. By [S2, Theorem 3.4],
$$\sum_{k=0}^{(p-1)/2}\f{\b{2k}k^3}{(-512)^k}\e
(-1)^{\f{p-1}4}(4x^2-2p)=(-1)^{\f{p-1}4}(-16y^2+2p)\mod {p^2}.\tag
2.14$$ By [L, Theorem 1.5],
$$\sum_{k=0}^{(p-1)/2}(6k+1)\f{\b{2k}k^3}{(-512)^k}\e
(-1)^{\f{p-1}4}p\mod {p^2}.$$ Thus,
$$\sum_{k=0}^{\f{p-1}2}\f{k\b{2k}k^3}{(-512)^k}
\e \f
16\Big((-1)^{\f{p-1}4}p-\sum_{k=0}^{\f{p-1}2}\f{\b{2k}k^3}{(-512)^k}\Big)
\e \f 16(-1)^{\f{p-1}4}(16y^2-p)\mod {p^2}.\tag 2.15$$ Set
$t=\f{3\sqrt 2-4}8$. Then $t(t+1)=\f 1{32}$ and $(2t+1)^2=\f 98$.
Now, from the above and Corollary 2.5 we deduce that
$$\align\sum_{k=0}^{(p-1)/2}\f{\b{2k}k^3}{(-512)^k(k+1)}
&\e \sum_{k=0}^{(p-1)/2}k\f{\b{2k}k^3}{(-512)^k}+\f{9/8}{1/32}\cdot
\f{(\sum_{k=0}^{(p-1)/2}k\b{2k}k^3(-512)^{-k})^2}
{\sum_{k=0}^{(p-1)/2}\b{2k}k^3(-512)^{-k}}
\\&\e (-1)^{\f{p-1}4}(-16y^2+2p)+\f{(16y^2-p)^2}{(-1)^{\f{p-1}4}(-16y^2+2p)}
\\&\e (-1)^{\f{p-1}4}(-32y^2+2p)\mod {p^2}.\endalign$$
This proves the theorem. \vskip 0.25cm
\pro{Lemma 2.3} Let $p$ be an
odd prime. If $p\e 1\mod 4$ and so $p=x^2+4y^2$ with $4\mid x-1$,
then
$$\sum_{k=0}^{(p-1)/2}\f{k\b{2k}k^2}{(-16)^k}
\e (-1)^{\f{p-1}4}\Big(-\f x2+\f p{4x}\Big)\mod {p^2}.$$ If $p\e
3\mod 4$, then
$$\align\sum_{k=0}^{(p-1)/2}\f{k\b{2k}k^2}{(-16)^k}
&\e\f
14(-1)^{\f{p+1}4}\Big(\b{(p-1)/2}{(p-3)/4}+\Big(\b{(p-1)/2}{(p-3)/4}
+\b{(p-1)/2}{(p-3)/4}^{-1}\\&\q-\f
12\cdot\f{2^{p-1}-1}p\b{(p-1)/2}{(p-3)/4}\Big)p\Big)\mod
{p^2}.\endalign$$
\endpro
Proof. It is clear that
$$\align &\sum_{k=0}^{p-1}(2k+1)\f{\b{2k}k^2}{(-16)^k}+\sum_{k=0}^{p-1}
\f{\b{2k}k^2}{(-16)^k(2k-1)}
\\&=\sum_{k=0}^{p-1}\f{4k^2\b{2k}k^2}{(-16)^k(2k-1)}
=\sum_{k=0}^{p-2}\f{4(k+1)^2\b{2(k+1)}{k+1}^2}{(-16)^{k+1}(2(k+1)-1)}
=-\sum_{k=0}^{p-2}(2k+1)\f{\b{2k}k^2}{(-16)^k}.\endalign$$ Thus,
$$\sum_{k=0}^{p-1}\f{k\b{2k}k^2}{(-16)^k}\e -\f
12\sum_{k=0}^{p-1}\f{\b{2k}k^2}{(-16)^k}-\f
14\sum_{k=0}^{p-1}\f{\b{2k}k^2}{(-16)^k(2k-1)}\mod {p^2}.$$ We first
assume that $p\e 1\mod 4$ and so $p=x^2+4y^2$ with $x\e 1\mod 4$.
 By [S1, Theorem 2.9],
$$\sum_{k=0}^{p-1}\f{\b{2k}k^2}{(-16)^k}\e
(-1)^{\f{p-1}4}\Big(2x-\f p{2x}\Big)\mod{p^2}.\tag 2.16$$ By [Su5,
Theorem 1.1], $$\sum_{k=0}^{p-1}\f{\b{2k}k^2}{(-16)^k(2k-1)}\e
-2(-1)^{\f{p-1}4}x\mod {p^2}.$$ Thus,
$$\sum_{k=0}^{p-1}\f{k\b{2k}k^2}{(-16)^k}
\e -\f 12(-1)^{\f{p-1}4}\Big(2x-\f p{2x}\Big)+\f 24(-1)^{\f{p-1}4}x
=(-1)^{\f{p-1}4}\Big(-\f x2+\f p{4x}\Big)\mod {p^2}.$$ This proves
the result in the case $p\e 1\mod 4$. Now assume that $p\e 3\mod 4$.
By [S11, Theorem 3.2],
$$\sum_{k=0}^{p-1}\f{\b{2k}k^2}{(-16)^k(2k-1)}
\e\f
12(-1)^{\f{p-3}4}(2p+3-2^{p-1})\b{\f{p-1}2}{\f{p-3}4}-(-1)^{\f{p-3}4}
\f p{\b{(p-1)/2}{(p-3)/4}}\mod {p^2}.$$ By [Su4, (1.6)],
$$\sum_{k=0}^{p-1}\f{\b{2k}k^2}{(-16)^k}\e (-1)^{\f{p-3}4}\f p{\b{(p-1)/2}{(p-3)/4}}\mod
{p^2}.\tag 2.17$$ Thus,
$$\align \sum_{k=0}^{p-1}\f{k\b{2k}k^2}{(-16)^k}
&\e -\f 12(-1)^{\f{p-3}4}\f p{\b{(p-1)/2}{(p-3)/4}} -\f 14\Big(\f
12(-1)^{\f{p-3}4}(2p+3-2^{p-1})\b{\f{p-1}2}{\f{p-3}4}
\\&\q-(-1)^{\f{p-3}4}
\f p{\b{(p-1)/2}{(p-3)/4}}\Big)\mod {p^2}. \endalign$$ This yields
the result in the case $p\e 3\mod 4$. The proof is now complete.
\vskip 0.25cm
 \pro{Theorem 2.7} Let $p$ be an odd prime. Then
$$\sum_{k=0}^{(p-1)/2}\f{\b{2k}k^3}{(-8)^k(k+1)}\e\cases
-24y^2+2p\mod {p^2}& \t{if $p=x^2+4y^2\e 1\mod 4$,}
\\\f 12R_1(p)+p\mod {p^2}&\t{if $p\e 3\mod 4$.}
\endcases$$
\endpro
Proof. Taking $t=1$ in Corollary 2.1 and then applying (2.16),
(2.17) and Lemma 2.3 gives
$$\align&\sum_{k=0}^{(p-1)/2}\f{\b{2k}k^3}{(-8)^k(k+1)}
\\&\e
\Big(\sum_{k=0}^{(p-1)/2}
\f{\b{2k}k^2}{(-16)^k}\Big)^2+8\Big(\sum_{k=0}^{(p-1)/2}
\f{k\b{2k}k^2}{(-16)^k}\Big)^2
\\&\e\cases 4x^2-2p+8
(-\f x2+\f p{4x})^2=-24y^2+2p\mod {p^2}\q\t{if $p=x^2+4y^2\e 1\mod
4$,}
\\0+\f
8{16}\Big(\b{(p-1)/2}{(p-3)/4}+\Big(\b{(p-1)/2}{(p-3)/4}
+\b{(p-1)/2}{(p-3)/4}^{-1}\\\ -\f
12\cdot\f{2^{p-1}-1}p\b{(p-1)/2}{(p-3)/4}\Big)p\Big)^2 \e
p+(p+1-2^{p-2})\b{(p-1)/2}{(p-3)/4}^2\mod
{p^2}\\\qq\qq\qq\qq\qq\qq\qq\qq\qq\qq\qq\ \t{if $4\mid p-3$.}
\endcases\endalign$$
This proves the theorem. \vskip 0.25cm

\par{\bf Remark 2.2} Let $p>3$ be a prime. In [S9], the author
conjectured that
$$\align&\sum_{k=0}^{(p-1)/2}\f{\b{2k}k^3}{(-8)^k(k+1)} \e
-24y^2+2p\mod {p^3}\q \t{for}\q p=x^2+4y^2\e 1\mod 4,
\\&(-1)^{[\f p4]}\sum_{k=0}^{\f{p-1}2}\f{\b{2k}k^3}{(-512)^k(k+1)}
\e\cases -32y^2+2p\mod {p^3}&\t{if $p=x^2+4y^2\e 1\mod 4$,}
\\-4R_1(p)-2p\mod {p^2}&\t{if $p\e 3\mod 4$,}
\endcases
\\&\sum_{k=0}^{(p-1)/2}\f{\b{2k}k^3}{16^k(k+1)}
\e\cases -16y^2+2p\mod {p^3}&\t{if $p=x^2+3y^2\e 1\mod 3$,}
\\-\f 43R_3(p)-\f 23p\mod {p^2}&\t{if $p\e 2\mod 3$,}
\endcases
\\&(-1)^{\f{p-1}2}\sum_{k=0}^{(p-1)/2}\f{\b{2k}k^3}{256^k(k+1)}
\e\cases -8y^2+2p\mod {p^3}&\t{if $p=x^2+3y^2\e 1\mod 3$,}
\\\f {16}3R_3(p)+\f 23p\mod {p^2}&\t{if $p\e 2\mod 3$,}
\endcases
\\&(-1)^{\f{p-1}2}\sum_{k=0}^{(p-1)/2}\f{\b{2k}k^3}{(-64)^k(k+1)}
\e -12y^2+2p\mod {p^3}\ \t{for}\ p=x^2+2y^2\e 1,3\mod 8,
\\&(-1)^{\f{p-1}2}\sum_{k=0}^{(p-1)/2}\f{\b{2k}k^3}{4096^k(k+1)}
\e 72y^2+2p\mod {p^3}\ \t{for}\ p=x^2+7y^2\e 1,2,4\mod 7.
\endalign$$

\pro{Lemma 2.4 ([WS, Theorems 1.2 and 1.3]} Let $p>3$ be a prime.
Then
$$\align&\sum_{k=0}^{p-1}\f{\b{2k}k\b{3k}k}{24^k}
\e\cases \b{(2p-2)/3}{(p-1)/3}\mod {p^2}&\t{if $3\mid p-1$,}
\\ p\b{(2p+2)/3}{(p+1)/3}^{-1}\mod {p^2}&\t{if $3\mid p-2$,}
\endcases
\\&\sum_{k=0}^{p-1}\f{k\b{2k}k\b{3k}k}{24^k}
\e\cases p\b{(2p-2)/3}{(p-1)/3}^{-1}-\b{(2p-2)/3}{(p-1)/3}\mod {p^2}
&\t{if $3\mid p-1$,}
\\-(p+1)\b{(2p+2)/3}{(p+1)/3}-p\b{(2p+2)/3}{(p+1)/3}^{-1}\mod {p^2}
&\t{if $3\mid p-2$.}
\endcases\endalign$$
\endpro
\pro{Theorem 2.8} Let $p>5$ be a prime. Then
$$\align
&\sum_{k=0}^{p-1}\f{k\b{2k}k^2\b{3k}k}{(-192)^k} \e\cases \f 35p-\f
15x^2\mod{p^2}&\t{if $3\mid p-1$ and so $4p=x^2+27y^2$,}
\\-\f p5\mod {p^2}&\t{if $p\e 2\mod 3$,}
\endcases
\\ &\sum_{k=0}^{p-2}\f{\b{2k}k^2\b{3k}k}{(-192)^k(k+1)} \e\cases \f
32x^2-4p\mod{p^2}\q\t{if $3\mid p-1$ and so $4p=x^2+27y^2$,}
\\2(2p+1)\b{[2p/3]}{[p/3]}^2+p\mod {p^2}\qq\t{if $p\e 2\mod 3$.}
\endcases\endalign$$
\endpro
Proof. For $p\e 1\mod 3$ it is known (see [BEW]) that
$$\b{2(p-1)/3}{(p-1)/3}\e \Ls x3\Big(\f px-x\Big)\mod {p^2}
\ \t{and so}\ \b{2(p-1)/3}{(p-1)/3}^2\e x^2-2p\mod {p^2}.$$ Taking
$a=-\f 13$ and $t=-\f 98$ in (2.2) and then applying (1.1) and Lemma
2.4 gives
$$\align\sum_{k=0}^{p-1}\f{k\b{2k}k^2\b{3k}k}{(-192)^k} &\e
\f 15\Big(\sum_{k=0}^{p-1}\f{\b{2k}k\b{3k}k}{24^k}\Big)
\Big(\sum_{k=0}^{p-1}\f{k\b{2k}k\b{3k}k}{24^k}\Big)
\\&\e \cases\f 15\b{(2p-2)/3}{(p-1)/3}\big(p\b{(2p-2)/3}{(p-1)/3}^{-1}
-\b{(2p-2)/3}{(p-1)/3}\big)\\\q\e \f p5-\f
15\b{(2p-2)/3}{(p-1)/3}^2\e\f {3p-x^2}5\mod
{p^2}\\\qq\qq\qq\qq\q\qq\q\t{if $3\mid p-1$ and $4p=x^2+27y^2$,}
\\-\f p5\b{(2p+2)/3}{(p+1)/3}^{-1}\big((p+1)\b{(2p+2)/3}
{(p+1)/3}+p\b{(2p+2)/3}{(p+1)/3}^{-1}\big)
\\\q\e -\f p5\mod {p^2}\qq\qq\qq\qq\qq\q\qq\t{if $3\mid
p-2$.}
\endcases\endalign$$
Taking $t=-\f 98$ in Corollary 2.2 and then applying Lemma 2.4 gives
$$\align&\sum_{k=0}^{p-2}\f{\b{2k}k^2\b{3k}k}{(-192)^k(k+1)}
\e \Big(\sum_{k=0}^{p-1}\f{\b{2k}k\b{3k}k}{24^k}\Big)^2 +\f 12
\Big(\sum_{k=0}^{p-1}\f{k\b{2k}k\b{3k}k}{24^k}\Big)^2
\\&\e\cases \b{(2p-2)/3}{(p-1)/3}^2+\f 12\Big(\b{(2p-2)/3}{(p-1)/3}^2-2p\Big)
\e \f 32x^2-4p\mod {p^2}
\\\qq\qq\qq\qq\qq\qq\qq\qq\q\t{if $3\mid p-1$ and $4p=x^2+27y^2$,}
\\ \f 12\Big(\b{2(p+1)/3}{(p+1)/3}+p\Big(\b{2(p+1)/3}{(p+1)/3}
+\b{2(p+1)/3}{(p+1)/3}^{-1}\Big)\Big)^2
\\\q\e \f
12\Big((2p+1)\b{2(p+1)/3}{(p+1)/3}^2+2p\Big)
=2(2p+1)\b{[2p/3]}{[p/3]}^2+p \mod
{p^2}\\\qq\qq\qq\qq\qq\qq\qq\qq\q\t{if $3\mid p-2$.}
\endcases\endalign$$
This completes the proof. \vskip 0.25cm

\pro{Lemma 2.5 ([WS])} Let $p>3$ be a prime. Then
$$\align&\sum_{k=0}^{p-1}\f{\b{2k}k\b{4k}{2k}}{48^k}
\e\cases \sls x3(2x-\f p{2x})\mod {p^2} &\t{if $p=x^2+3y^2\e 1\mod
3$,}
\\\f{3p}{2\b{(p+1)/2}{(p+1)/6}}\mod {p^2}&\t{if $p\e 2\mod 3$,}
\endcases
\\&\sum_{k=0}^{p-1}\f{k\b{2k}k\b{4k}{2k}}{48^k}
\\&\e\cases \sls x3(-x+\f p{2x})\mod {p^2}
\qq\qq\q\qq\t{if $p=x^2+3y^2\e 1\mod 3$,}\\\Big(-\f 13-\f p3-\f
29(2^{p-1}-1)+\f 14(3^{p-1}-1)\Big)\b{(p+1)/2}{(p+1)/6}-\f
{3p}{4\b{(p+1)/2}{(p+1)/6}} \mod
{p^2}\\\qq\qq\qq\qq\qq\qq\qq\qq\q\;\t{if $p\e 2\mod 3$.}
\endcases
\endalign$$
\pro{Theorem 2.9} Let $p>3$ be a prime. Then
$$\align&\sum_{k=0}^{p-1}\f{k\b{2k}k^2\b{4k}{2k}}{(-144)^k} \e\cases
-\f 45x^2+\f 35p\mod {p^2}&\t{if $p=x^2+3y^2\e 1\mod 3$,}
\\-\f p5\mod {p^2}&\t{if $p\e 2\mod 3$,}
\endcases
\\&\sum_{k=0}^{p-2}\f{\b{2k}k^2\b{4k}{2k}}{(-144)^k(k+1)} \e\cases
-16y^2+2p\mod {p^2}&\t{if $p=x^2+3y^2\e 1\mod 3$,}
\\\f 43R_3(p)+\f 23p\mod {p^2}&\t{if $p\e 2\mod 3$.}
\endcases
\endalign$$
\endpro
Proof. Putting $a=-\f 14$ and $t=-\f 43$ in (2.2) and then applying
Lemma 2.5 gives
$$\align&\sum_{k=0}^{p-1}\f{k\b{2k}k^2\b{4k}{2k}}{(-144)^k}
\\&\e \f 25\Big(\sum_{k=0}^{p-1}\f{\b{2k}k\b{4k}{2k}}{48^k}\Big)
\Big(\sum_{k=0}^{p-1}\f{k\b{2k}k\b{4k}{2k}}{48^k}\Big)
\\&\e \cases\f 25(2x-\f p{2x})(-x+\f p{2x})\e -\f 45x^2+\f 35p
\mod {p^2}&\t{if $p=x^2+3y^2\e 1\mod 3$,}
\\\f{3p}{2\b{(p+1)/2}{(p+1)/6}}\Big(-\f 13\b{(p+1)/2}{(p+1)/6}\Big)
=-\f p2\mod {p^2}&\t{if $p\e 2\mod 3$,}
\endcases
\endalign$$
and taking $t=-\f 43$ in Corollary 2.3 and then applying Lemma 2.5
gives
$$\align&\sum_{k=0}^{p-2}\f{k\b{2k}k^2\b{4k}{2k}}{(-144)^k(k+1)}
\\&\e \Big(\sum_{k=0}^{p-1}\f{\b{2k}k\b{4k}{2k}}{48^k}\Big)^2
+\f 43\Big(\sum_{k=0}^{p-1}\f{k\b{2k}k\b{4k}{2k}}{48^k}\Big)^2
\\&\e\cases(2x-\f p{2x})^2+\f 43(-x+\f p{2x})^2\e -16y^2+2p\mod {p^2}
\q\t{if $p=x^2+3y^2\e 1\mod 3$,}
\\\f 43\Big(\f 19\b{(p+1)/2}{(p+1)/6}^2-\f 23p\Big(\big(-\f 13-\f 29\cdot
\f{2^{p-1}-1}p+\f
14\cdot\f{3^{p-1}-1}p\big)\b{(p+1)/2}{(p+1)/6}^2-\f 34\Big)\Big)
\\\q=\f 23p+\f 43R_3(p)\mod {p^2}\q\qq\qq\qq\qq\qq\ \t{if $p\e 2\mod 3$}.
\endcases\endalign$$
The proof is now complete.

\pro{Lemma 2.6 ([WS, Theorem 5.2])} Let $p>3$ be a prime. Then
$$\align&\sum_{k=0}^{p-1}\f{\b{2k}k\b{4k}{2k}}{72^k}
\e\cases \sls 6p\b{\f{p-1}2}{\f{p-1}4}\big(1-\f
12(2^{p-1}-1)\big)\mod {p^2} &\t{if $p=x^2+4y^2\e 1\mod 4$,}
\\\ls 6p\f{p}{3\b{(p-1)/2}{(p-3)/4}}\mod {p^2}&\t{if $p\e 3\mod 4$,}
\endcases
\\&\sum_{k=0}^{p-1}\f{k\b{2k}k\b{4k}{2k}}{72^k}
\e\cases -\f 12\ls 6p\f p{\b{(p-1)/2}{(p-1)/4}} +\f 12\ls
6p\b{(p-1)/2}{(p-1)/4}\big(1-\f 12(2^{p-1}-1)\big)\mod {p^2}
\\\qq\qq\qq\qq\qq\qq\q\qq\qq\t{if $p=x^2+4y^2\e 1\mod 4$,}\\\ls 6p\f
{p}{6\b{(p-1)/2}{(p-3)/4}}+\ls 6p \big(\f 32+\f 32p-\f
34(2^{p-1}-1)\big)\b{(p-1)/2}{(p-3)/4} \mod
{p^2}\\\qq\qq\qq\qq\qq\qq\qq\qq\q\;\t{if $p\e 3\mod 4$.}
\endcases
\endalign$$
\pro{Theorem 2.10} Let $p$ be a prime such that $p\not=2,3,7$. Then
$$\align
&\sum_{k=0}^{p-1}\f{k\b{2k}k^2\b{4k}{2k}}{648^k} \e\cases \f 37p-\f
47x^2\mod {p^2}&\t{if $p=x^2+4y^2\e 1\mod 4$,}
\\-\f p7\mod {p^2}&\t{if $p\e 3\mod 4$,}
\endcases
\\&\sum_{k=0}^{p-2}\f{\b{2k}k^2\b{4k}{2k}}{648^k(k+1)} \e\cases -\f
{40}3y^2+2p\mod {p^2}&\t{if $p=x^2+4y^2\e 1\mod 4$,}
\\-\f 32R_1(p)-\f p3\mod {p^2}&\t{if $p\e 3\mod 4$.}
\endcases
\endalign$$
\endpro
Proof. For $p\e 1\mod 4$ and so $p=x^2+4y^2$, from [S2, Lemma 3.4]
we know that
$$\b{(p-1)/2}{(p-1)/4}^2\e 2^{p-1}(4x^2-2p)\mod {p^2}.\tag 2.18$$
Now, taking $a=-\f 14$ and $t=-\f 89$ in (2.2) and then applying
Lemma 2.6 and (2.18) gives
$$\align
&\sum_{k=0}^{p-1}\f{k\b{2k}k^2\b{4k}{2k}}{648^k}
\\&\e -\f 27\Big(\sum_{k=0}^{p-1}\f{\b{2k}k\b{4k}{2k}}{72^k}\Big)
\Big(\sum_{k=0}^{p-1}\f{k\b{2k}k\b{4k}{2k}}{72^k}\Big)
\\&\e\cases
-\f 27\Big(1-\f 12(2^{p-1}-1)\Big)\Big(\f
12\b{(p-1)/2}{(p-1)/4}^2-\f 12p \Big(1+\f
12\cdot\f{2^{p-1}-1}p\b{(p-1)/2}{(p-1)/4}^2\Big)\Big)
\\\q\e-\f 27\Big(\f 12\b{(p-1)/2}{(p-1)/4}^2(1-(2^{p-1}-1))-\f p2\Big)
\e -\f 47x^2+\f 37p\mod {p^2}\\\qq\qq\qq\qq\qq\qq\qq\qq\qq \q\t{if
$p=x^2+4y^2\e 1\mod 4$,} \\ -\f 27\ls 6p\f
p{3\b{(p-1)/2}{(p-3)/4}}\cdot\f 32\ls 6p\b{(p-1)/2}{(p-3)/4}=-\f
p7\mod {p^2}\q\t{if $p\e 3\mod 4$,}
\endcases
\endalign$$
and taking $t=-\f 89$ in Corollary 2.3 and then applying Lemma 2.6
gives
$$\align &\sum_{k=0}^{p-2}\f{\b{2k}k^2\b{4k}{2k}}{648^k(k+1)}
\\&\e \Big(\sum_{k=0}^{p-1}\f{\b{2k}k\b{4k}{2k}}{72^k}\Big)^2-\f 23
\Big(\sum_{k=0}^{p-1}\f{k\b{2k}k\b{4k}{2k}}{72^k}\Big)^2
\\&\e\cases \b{\f{p-1}2}{\f{p-1}4}^2(1-(2^{p-1}-1))-\f 23
\Big(\f 14\b{\f{p-1}2}{\f{p-1}4}^2-\f p2\Big(1+\f
{2^{p-1}-1}{2p}\b{\f{p-1}2}{\f{p-1}4}^2\Big)\Big)
\\\q\e \f p3+\f 56\b{\f{p-1}2}{\f{p-1}4}^2(1-(2^{p-1}-1))
\e \f p3+\f 56(4x^2-2p)=-\f{40}3y^2+2p\mod
{p^2}\\\qq\qq\qq\qq\qq\qq\q\qq\qq\qq\q\t{if $p=x^2+4y^2 \e 1\mod
4$,}
\\0-\f 23\Big(\f 32\b{(p-1)/2}{(p-3)/4}+p\Big(\f
1{6\b{(p-1)/2}{(p-3)/4}}+\big(\f 32-\f
34\cdot\f{2^{p-1}-1}p\big)\b{(p-1)/2}{(p-3)/4}\Big)\Big)^2
\\\q\e -\f p3-\f 32R_1(p)\mod {p^2}\qq\q\q\qq\qq\t{if $p\e 3\mod 4$.}
\endcases
\endalign$$
This proves the theorem. \vskip 0.25cm
\par{\bf Remark 2.3} Let $p>3$ be a prime.
By [S2, Theorem 5.1], [S3, Theorem 4.2] and [WS, Corollaries 1.1,
5.1 and 5.2],
$$\align&\sum_{k=0}^{p-1}\f{\b{2k}k^2\b{3k}k}{(-192)^k}\e\cases x^2-2p\mod
{p^2}&\t{if $3\mid p-1$ and $4p=x^2+27y^2$,}
\\0\mod {p^2}&\t{if $3\mid p-2$,}
\endcases\tag 2.19
\\&\sum_{k=0}^{p-1}\f{\b{2k}k^2\b{4k}{2k}}{(-144)^k} \e\cases
4x^2-2p\mod {p^2}&\t{if $p=x^2+3y^2\e 1\mod 3$,}
\\0\mod {p^2}&\t{if $p\e 2\mod 3$,}\endcases\tag 2.20
\\&\sum_{k=0}^{p-1}\f{\b{2k}k^2\b{4k}{2k}}{648^k} \e
\cases 4x^2-2p\mod {p^2}&\t{if $p=x^2+4y^2\e 1\mod 4$,}
\\0\mod {p^2}&\t{if $p\e 3\mod 4$.}\endcases\tag 2.21\endalign$$
In [S7] the author conjectured that
$$\align&\sum_{k=0}^{p-1}\f{\b{2k}k^2\b{3k}k}{(-192)^k}\e\cases x^2-2p
-\f{p^2}{x^2}\mod {p^3}&\t{if $3\mid p-1$ and $4p=x^2+27y^2$,}
\\\f 34p^2\b{(2p-1)/3}{(p-2)/3}^{-2}\mod {p^3}&\t{if $3\mid p-2$,}
\endcases
\\&\sum_{k=0}^{p-1}\f{\b{2k}k^2\b{4k}{2k}}{(-144)^k} \e\cases
4x^2-2p-\f{p^2}{4x^2}\mod {p^3}&\t{if $p=x^2+3y^2\e 1\mod 3$,}
\\p^2\b{(p-1)/2}{(p-5)/6}^{-2}\mod {p^3}&\t{if $p\e 2\mod 3$,}\endcases
\\&\sum_{k=0}^{p-1}\f{\b{2k}k^2\b{4k}{2k}}{648^k} \e
\cases 4x^2-2p-\f{p^2}{4x^2}\mod {p^3}&\t{if $p=x^2+4y^2\e 1\mod
4$,}
\\-\f 5{36}p^2\b{(p-3)/2}{(p-3)/4}^{-2}
\mod {p^3}&\t{if $p\e 3\mod 4$.}\endcases\endalign$$
 In [S9], the author conjectured that
$$\align&\sum_{k=0}^{p-2}\f{\b{2k}k^2\b{3k}k}{(-192)^k(k+1)}\e \f 32x^2-4p
\mod {p^3}\q\t{for $3\mid p-1$ and so $4p=x^2+27y^2$,}
\\&\sum_{k=0}^{p-2}\f{\b{2k}k^2\b{4k}{2k}}{(-144)^k(k+1)} \e
-16y^2+2p\mod {p^3}\q\t{for $p=x^2+3y^2\e 1\mod 3$,}
\\&\sum_{k=0}^{p-2}\f{\b{2k}k^2\b{4k}{2k}}{648^k(k+1)} \e -\f
{40}3y^2+2p\mod {p^3}\q\t{for $p=x^2+4y^2\e 1\mod 4$.}\endalign$$ We
also conjecture that
$$\align&\sum_{k=0}^{p-1}(5k+1)\f{\b{2k}k^2\b{4k}{2k}}{(-144)^k}
\e \Ls p3p+\f 52p^3U_{p-3}\mod {p^4},\tag 2.22
\\&\sum_{k=0}^{p-1}(7k+1)\f{\b{2k}k^2\b{4k}{2k}}{648^k}
\e (-1)^{\f{p-1}2}p-\f{745}{447}p^3E_{p-3}\mod {p^4}\q\t{for
$p\not=149$},\tag 2.23\endalign$$ where $\{E_n\}$ (Euler numbers)
and $\{U_n\}$ are given by
$$\align &E_{2n-1}=0,\ E_0=1,\ E_{2n}=-\sum_{k=1}^n\b {2n}{2k}E_{2n-2k}\
(n\ge 1),\\& U_{2n-1}=0,\ U_0=1,\ U_{2n}=-2\sum_{k=1}^n\b
{2n}{2k}U_{2n-2k}\ (n\ge 1).\endalign$$ In [Su1], Z.W. Sun made a
conjecture equivalent to
$$\sum_{k=0}^{p-1}(5k+1)\f{\b{2k}k^2\b{3k}k}{(-192)^k}\e\Ls p3p+\f 53p^3U_{p-3}
\mod {p^4}.$$

\section*{3. Congruences for
$\sum_{k=0}^{p-1}k^2\b ak\b{-1-a}k\b{2k}k\f 1{m^k}$ and
\\$\sum_{k=0}^{p-1}k^3\b ak\b{-1-a}k\b{2k}k\f 1{m^k}$
 modulo $p^2$}
\pro{Theorem 3.1} Let $p$ be an odd prime, $a,m\in\Bbb Z_p$,
$a\not\e 0,-1\mod p$ and $m\not\e 0\mod p$. Then
$$\align \f{m-4}2&\sum_{k=0}^{p-1}\f{k^2\b ak\b{-1-a}k\b{2k}k}{m^k}
\e \sum_{k=0}^{p-1}\f{k\b ak\b{-1-a}k\b{2k}k}{m^k}
-2a(a+1)\sum_{k=0}^{p-1}\f{\b ak\b{-1-a}k\b{2k}k}{m^k}
\\&\qq\qq\qq\qq\qq+a(a+1)\sum_{k=0}^{p-2} \f{\b ak\b{-1-a}k\b{2k}k}{m^k(k+1)}\mod {p^3},
\\\f{m-4}2&\sum_{k=0}^{p-1}\f{k^3\b ak\b{-1-a}k\b{2k}k}{m^k}
\e 3\sum_{k=0}^{p-1}\f{k^2\b ak\b{-1-a}k\b{2k}k}{m^k}
-(2a(a+1)-1)\sum_{k=0}^{p-1}\f{k\b ak\b{-1-a}k\b{2k}k}{m^k}
\\&\qq\qq\qq\qq\qq-a(a+1)\sum_{k=0}^{p-1} \f{\b ak\b{-1-a}k\b{2k}k}{m^k}\mod {p^3}.
\endalign$$
Thus, for $m\not\e 4\mod p$,
$$\align&\sum_{k=0}^{p-1}\f{k^3\b ak\b{-1-a}k\b{2k}k}{m^k}
\\&\e \f{(2-4a(a+1))(m-4)+12}{(m-4)^2}\sum_{k=0}^{p-1}\f{k\b ak\b{-1-a}k\b{2k}k}{m^k}
-\f{2a(a+1)(m+8)}{(m-4)^2}\sum_{k=0}^{p-1}\f{\b
ak\b{-1-a}k\b{2k}k}{m^k}\\&\qq+\f{12a(a+1)}{(m-4)^2}
\sum_{k=0}^{p-2}\f{\b ak\b{-1-a}k\b{2k}k}{m^k(k+1)}\mod{p^3}.
\endalign$$

\endpro
Proof. Since
$$\aligned &(k+1)^2\b a{k+1}\b{-1-a}{k+1}\b{2(k+1)}{k+1}
\\&=(k+1)^2\cdot\f{a-k}{k+1}\b ak\cdot\f{-1-a-k}{k+1}\b{-1-a}k\cdot
\f{2(2k+1)}{k+1}\b{2k}k
\\&=\Big(4k^2+2k-4a(a+1)+\f{2a(a+1)}{k+1}\Big)\b ak\b{-1-a}k\b{2k}k,
\endaligned\tag 3.1$$
we see that
$$\align\sum_{k=0}^{p-1}\f{k^2\b ak\b{-1-a}k\b{2k}k}{m^k}
&=\sum_{k=0}^{p-2}\f{(k+1)^2\b
a{k+1}\b{-1-a}{k+1}\b{2(k+1)}{k+1}}{m^{k+1}} \\&=\f 1m
\sum_{k=0}^{p-2} \Big(4k^2+2k-4a(a+1)+\f{2a(a+1)}{k+1}\Big)\f{\b
ak\b{-1-a}k\b{2k}k}{m^k}
\endalign$$
and
$$\align\sum_{k=0}^{p-1}\f{k^3\b ak\b{-1-a}k\b{2k}k}{m^k}
&=\sum_{k=0}^{p-2}\f{(k+1)^3\b
a{k+1}\b{-1-a}{k+1}\b{2(k+1)}{k+1}}{m^{k+1}} \\&=\f 1m
\sum_{k=0}^{p-2} \big((4k^2+2k-4a(a+1))(k+1)+2a(a+1)\big)\f{\b
ak\b{-1-a}k\b{2k}k}{m^k}
\\&=\f 1m
\sum_{k=0}^{p-2} \big(4k^3+6k^2+(2-4a(a+1))k-2a(a+1)\big)\f{\b
ak\b{-1-a}k\b{2k}k}{m^k}.
\endalign$$
Note that $\b a{p-1}\e \b{-1-a}{p-1}\e \b{2p-2}{p-1}\e 0\mod p$.
From the above we obtain the first two congruences in Theorem 3.1.
Combining the two congruences we deduce the remaining result.
\par\q
\pro{Corollary 3.1} Let $p$ be an odd  prime, $m\in\Bbb Z_p$ and
$m\not\e 0\mod p$. Then
$$\f{m-4}2\sum_{k=0}^{\f{p-1}2}\f{k^2\b{2k}k^3}{(16m)^k}
\e \sum_{k=0}^{\f{p-1}2}\f{k\b{2k}k^3}{(16m)^k}+\f 12
\sum_{k=0}^{\f{p-1}2}\f{\b{2k}k^3}{(16m)^k}-\f
14\sum_{k=0}^{\f{p-1}2}\f{\b{2k}k^3}{(16m)^k(k+1)}\mod {p^3}$$ and
for $m\not\e 4\mod p$,
$$\align \sum_{k=0}^{\f{p-1}2}\f{k^3\b{2k}k^3}{(16m)^k} &\e
\f{3m}{(m-4)^2}\sum_{k=0}^{\f{p-1}2}\f{k\b{2k}k^3}{(16m)^k}+\f
{m+8}{2(m-4)^2} \sum_{k=0}^{\f{p-1}2}\f{\b{2k}k^3}{(16m)^k}\\&\qq-\f
3{(m-4)^2}\sum_{k=0}^{\f{p-1}2}\f{\b{2k}k^3}{(16m)^k(k+1)}\mod
{p^3}.
\endalign$$
\endpro
Proof. Note that $p\mid\b{2k}k$ for $\f p2<k<p$ and
$\b{-1/2}k=\b{2k}k(-4)^{-k}$. Taking $a=-\f 12$ in Theorem 3.1
yields the result. \vskip 0.25cm
 \pro{Theorem 3.2} Let $p$ be an odd
prime such that $p\e 1,2,4\mod 7$ and so $p=x^2+7y^2$. Then
$$\align&\sum_{k=0}^{(p-1)/2}k^2\b{2k}k^3\e
 \f{736x^2}{1323}-\f{272p}{441}\mod {p^2},
\\&\sum_{k=0}^{(p-1)/2}k^3\b{2k}k^3\e
-\f{5408x^2}{27783}+\f{2992p}{9261}\mod {p^2},
\\&(-1)^{\f{p-1}2}\sum_{k=0}^{(p-1)/2}k^2\f{\b{2k}k^3}{4096^k}
\e  \f{43x^2}{1323}-\f{13}{441}p\mod {p^2},
\\&(-1)^{\f{p-1}2}\sum_{k=0}^{(p-1)/2}k^3\f{\b{2k}k^3}{4096^k}
\e\f{169x^2}{55566}-\f{31}{74088}p\mod {p^2}.
\endalign$$
\endpro
Proof. Taking $m=\f 1{16}$ in Corollary 3.1 and then applying (2.4),
(2.5) and Theorem 2.3 we deduce that
$$\align\sum_{k=0}^{(p-1)/2}k^2\b{2k}k^3&\e-\f{32}{63}\Big(
\sum_{k=0}^{(p-1)/2}k\b{2k}k^3+\f 12\sum_{k=0}^{(p-1)/2}\b{2k}k^3-\f
14\sum_{k=0}^{(p-1)/2}\f{\b{2k}k^3}{k+1}\Big)
\\&\e-\f {32}{63}\Big(\f 8{21}(3p-4x^2)+\f 12(4x^2-2p)-\f
14\Big(-44\f{p-x^2}7+2p\Big)\Big) \\&=
 \f{736x^2}{1323}-\f{272p}{441}\mod {p^2}.
 \endalign$$
 and
$$\align\sum_{k=0}^{(p-1)/2}k^3\b{2k}k^3&\e\f{16^2}{63^2}\Big(\f 3{16}
\sum_{k=0}^{(p-1)/2}k\b{2k}k^3+\f
{129}{32}\sum_{k=0}^{(p-1)/2}\b{2k}k^3-3\sum_{k=0}^{(p-1)/2}\f{\b{2k}k^3}{k+1}\Big)
\\&\e\f {16^2}{63^2}\Big(\f 3{16}\cdot
\f 8{21}(3p-4x^2)+\f
{129}{32}(4x^2-2p)-3\Big(-44\f{p-x^2}7+2p\Big)\Big) \\&=
-\f{5408x^2}{27783}+\f{2992p}{9261}\mod {p^2}.
 \endalign$$
Similarly, taking $m=256$ in Corollary 3.1 and then applying (2.6),
(2.7) and Theorem 2.3 we deduce the remaining congruences. \vskip
0.25cm
\par{\bf Remark 3.1} Let $p$ be an odd prime such that $p\e 1,2,4\mod 7$
and so $p=x^2+7y^2$. In [S8], the author conjectured that
$$\align&\sum_{k=0}^{p-1}k^2\b{2k}k^3\e
\f{736x^2}{1323}-\f{272p}{441}+\f{20p^2}{1323x^2}\mod {p^3},
\\&\sum_{k=0}^{p-1}k^3\b{2k}k^3\e
-\f{5408x^2}{27783}+\f{2992p}{9261}-\f{1774p^2}{27783x^2}\mod {p^3},
\\&(-1)^{\f{p-1}2}\sum_{k=0}^{(p-1)/2}k^2\f{\b{2k}k^3}{4096^k}
\e \f{43x^2}{1323}-\f{13}{441}p-\f{p^2}{1323x^2}\mod {p^3},
\\&(-1)^{\f{p-1}2}\sum_{k=0}^{(p-1)/2}k^3\f{\b{2k}k^3}{4096^k}
\e \f{169x^2}{55566}-\f{31}{74088}p-\f{71p^2}{444528x^2}\mod
{p^3}.\endalign$$
 \vskip 0.25cm
 \pro{Theorem 3.3} Let $p$ be a prime such that $p\e 1\mod 3$ and so
 $p=x^2+3y^2$. Then
 $$\align &\sum_{k=0}^{(p-1)/2}\f{k^2\b{2k}k^3}{16^k}\e
 \f 49x^2-\f 59 p\mod {p^2},
 \\&\sum_{k=0}^{(p-1)/2}\f{k^3\b{2k}k^3}{16^k}\e
-\f 29x^2+\f 49 p\mod {p^2},
\\&(-1)^{\f{p-1}2}\sum_{k=0}^{(p-1)/2}\f{k^2\b{2k}k^3}{256^k} \e \f
19x^2-\f p{18}\mod {p^2},
\\&(-1)^{\f{p-1}2}\sum_{k=0}^{(p-1)/2}\f{k^3\b{2k}k^3}{256^k}
\e \f 1{18}x^2+\f p{72}\mod{p^2}.\endalign$$
\endpro
Proof. Taking $m=1,16$ in Corollary 3.1 and then applying
(2.8)-(2.11) and Theorem 2.4 yields the result. \vskip 0.25cm
 \pro{Theorem 3.4} Let $p$ be a prime such that $p\e
1,3\mod 8$ and so
 $p=x^2+2y^2$. Then
 $$\align &(-1)^{\f{p-1}2}\sum_{k=0}^{(p-1)/2}\f{k^2\b{2k}k^3}{(-64)^k}\e
 \f 18x^2-\f 3{16}p\mod {p^2},
\\&(-1)^{\f{p-1}2}\sum_{k=0}^{(p-1)/2}\f{k^3\b{2k}k^3}{(-64)^k} \e
\f 1{32}x^2-\f p{64}\mod {p^2}.\endalign$$
\endpro
Proof. Taking $m=-4$ in Corollary 3.1 and then applying
(2.12)-(2.13) and Theorem 2.5 yields the result. \vskip 0.25cm
\pro{Theorem 3.5} Let $p$ be a prime of the form $4k+1$ and so
$p=x^2+4y^2$. Then
$$\align&(-1)^{\f{p-1}4}\sum_{k=0}^{(p-1)/2}\f{k^2\b{2k}k^3}{(-512)^k} \e \f
{x^2}{27}-\f p{18}\mod {p^2},
\\&(-1)^{\f {p-1}4}\sum_{k=0}^{(p-1)/2}\f{k^3\b{2k}k^3}{(-512)^k}
\e -\f {x^2}{162}-\f p{216}\mod {p^2}.
\endalign$$
\endpro
Proof. Putting $m=-32$ in Corollary 3.1 and then applying
(2.14)-(2.15) and Theorem 2.6 yields the result. \vskip 0.25cm
 \pro{Theorem 3.6} Let $p>3$ be a prime. Then
$$\align&\sum_{k=0}^{\f{p-1}2}\f{k^2\b{2k}k^3}{(-8)^k}\e\cases
\f{10}{27}x^2-\f 49p\mod {p^2}& \t{if $p=x^2+4y^2\e 1\mod 4$,}
\\\f 7{27}p+\f 1{18}R_1(p)\mod {p^2}&\t{if $p\e 3\mod 4$,}
\endcases
\\&\sum_{k=0}^{\f{p-1}2}\f{k^3\b{2k}k^3}{(-8)^k}\e\cases
-\f{4}{81}x^2+\f 4{27}p\mod {p^2}& \t{if $p=x^2+4y^2\e 1\mod 4$,}
\\-\f {10}{81}p-\f 2{27}R_1(p)\mod {p^2}&\t{if $p\e 3\mod 4$.}
\endcases
\endalign$$
\endpro
Proof. By [S2, Theorems 3.3 and 3.4],
$$\sum_{k=0}^{(p-1)/2}\f{\b{2k}k^3}{(-8)^k}\e
\cases 4x^2-2p \mod {p^2}&\t{if $p=x^2+4y^2\e 1\mod 4$,}
\\0\mod {p^2}&\t{if $p\e 3\mod 4$.}\endcases\tag
3.2$$ By [GZ],
$$\sum_{k=0}^{(p-1)/2}(3k+1)\f{\b{2k}k^3}{(-8)^k}\e
(-1)^{\f{p-1}2}p\mod {p^3}.$$ This together with $(3.2)$ gives
$$\sum_{k=0}^{\f{p-1}2}\f{k\b{2k}k^3}{(-8)^k}\e
\cases p-\f 43x^2 \mod {p^2}&\t{if $p=x^2+4y^2\e 1\mod 4$,}
\\-\f 13p\mod {p^2}&\t{if $p\e 3\mod 4$.}\endcases\tag
3.3$$ Taking $m=-\f 12$ in Corollary 3.1 gives
$$\align&\sum_{k=0}^{\f{p-1}2}\f{k^2\b{2k}k^3}{(-8)^k}
\e -\f 49\sum_{k=0}^{\f{p-1}2}\f{k\b{2k}k^3}{(-8)^k} -\f
29\sum_{k=0}^{\f{p-1}2}\f{\b{2k}k^3}{(-8)^k}+\f
19\sum_{k=0}^{\f{p-1}2}\f{\b{2k}k^3}{(-8)^k(k+1)}\mod{p^2},
\\&\sum_{k=0}^{\f{p-1}2}\f{k^3\b{2k}k^3}{(-8)^k}
\e -\f 2{27}\sum_{k=0}^{\f{p-1}2}\f{k\b{2k}k^3}{(-8)^k} +\f
5{27}\sum_{k=0}^{\f{p-1}2}\f{\b{2k}k^3}{(-8)^k}-\f
4{27}\sum_{k=0}^{\f{p-1}2}\f{\b{2k}k^3}{(-8)^k(k+1)}
\mod{p^2}.\endalign$$ Now applying (3.2), (3.3) and Theorem 2.7
yields the result. \vskip 0.25cm

 \pro{Theorem 3.7} Let $p>5$ be a
prime. Then
$$\align
&\sum_{k=0}^{p-1}\f{k^2\b{2k}k^2\b{3k}{k}}{(-192)^k} \e\cases \f
{2}{125}x^2-\f {27p}{250}\mod{p^2}\q\t{if $3\mid p-1$ and so
$4p=x^2+27y^2$,}
\\\f 2{25}(2p+1)\b{[2p/3]}{[p/3]}^2+\f {19}{250}p\mod {p^2}
\q\t{if $p\e 2\mod{3}$,}
\endcases
\\&\sum_{k=0}^{p-1}\f{k^3\b{2k}k^2\b{3k}{k}}{(-192)^k}
\e\cases \f {21}{6250}x^2-\f {221p}{12500}\mod{p^2}\q \t{if $3\mid
p-1$ and so $4p=x^2+27y^2$,}
\\-\f {27}{625}(2p+1)\b{[2p/3]}{[p/3]}^2+\f {137}{12500}p\mod {p^2}
\q\t{if $p\e 2\mod{3}$.}
\endcases
\endalign$$
\endpro
Proof. Putting $a=-\f 13$ and $m=-\f{64}9$ in Theorem 3.1 and then
applying (1.1) gives
$$\align& \sum_{k=0}^{p-1}\f{k^2\b{2k}k^2\b{3k}{k}}{(-192)^k}\\&\e
-\f{9}{50}\sum_{k=0}^{p-1}\f{k\b{2k}k^2\b{3k}{k}}{(-192)^k} -\f
2{25}\sum_{k=0}^{p-1}\f{\b{2k}k^2\b{3k}{k}}{(-192)^k} +\f
1{25}\sum_{k=0}^{p-2}\f{\b{2k}k^2\b{3k}{k}}{(-192)^k(k+1)}\mod{p^3},
\\& \sum_{k=0}^{p-1}\f{k^3\b{2k}k^2\b{3k}{k}}{(-192)^k}\\&\e
-\f{407}{2500}\sum_{k=0}^{p-1}\f{k\b{2k}k^2\b{3k}{k}}{(-192)^k} +\f
2{625}\sum_{k=0}^{p-1}\f{\b{2k}k^2\b{3k}{k}}{(-192)^k} -\f
{27}{1250}\sum_{k=0}^{p-2}\f{\b{2k}k^2\b{3k}{k}}{(-192)^k(k+1)}\mod{p^3}.
\endalign$$
Now applying Theorem 2.8 and (2.19) yields the result. \vskip 0.25cm

\pro{Theorem 3.8} Let $p>5$ be a prime. Then
$$\align&\sum_{k=0}^{p-1}\f{k^2\b{2k}k^2\b{4k}{2k}}{(-144)^k} \e\cases
\f {12}{125}x^2-\f {19}{125}p\mod {p^2}&\t{if $p=x^2+3y^2\e 1\mod
3$,}
\\\f 2{25}R_3(p)+\f {13}{125}p\mod {p^2}&\t{if $p\e 2\mod 3$,}
\endcases
\\&\sum_{k=0}^{p-1}\f{k^3\b{2k}k^2\b{4k}{2k}}{(-144)^k} \e\cases
\f {62}{3125}x^2+\f {6}{3125}p\mod {p^2}&\t{if $p=x^2+3y^2\e 1\mod
3$,}
\\-\f {48}{625}R_3(p)-\f {37}{3125}p\mod {p^2}&\t{if $p\e 2\mod 3$.}
\endcases
\endalign$$
\endpro
Proof. Taking $a=-\f 14$ and $m=-\f 94$ in Theorem 3.1 and then
applying (1.1) gives
$$\align &\sum_{k=0}^{p-1}\f{k^2\b{2k}k^2\b{4k}{2k}}{(-144)^k}
\\&\e -\f 8{25}\sum_{k=0}^{p-1}\f{k\b{2k}k^2\b{4k}{2k}}{(-144)^k}-
\f 3{25}\sum_{k=0}^{p-1}\f{\b{2k}k^2\b{4k}{2k}}{(-144)^k}+\f
3{50}\sum_{k=0}^{p-2}\f{\b{2k}k^2\b{4k}{2k}}{(-144)^k(k+1)}
 \mod {p^3}\endalign$$
 and
$$\align &\sum_{k=0}^{p-1}\f{k^3\b{2k}k^2\b{4k}{2k}}{(-144)^k}
\\&\e -\f {83}{625}\sum_{k=0}^{p-1}\f{k\b{2k}k^2\b{4k}{2k}}{(-144)^k}+
\f{69}{1250}\sum_{k=0}^{p-1}\f{\b{2k}k^2\b{4k}{2k}}{(-144)^k}-\f
{36}{625}\sum_{k=0}^{p-2}\f{\b{2k}k^2\b{4k}{2k}}{(-144)^k(k+1)}
 \mod {p^3}.\endalign$$
Now applying Theorem 2.9 and (2.20) yields the result. \vskip 0.25cm
\pro{Theorem 3.9} Let $p$ be a prime with $p\not=2,3,7$. Then
$$\align
&\sum_{k=0}^{p-1}\f{k^2\b{2k}k^2\b{4k}{2k}}{648^k} \e \cases
\f{34}{343}x^2-\f{8p}{343}\mod {p^2}&\t{if $p=x^2+4y^2\e 1\mod 4$,}
\\\f {9}{98}R_1(p)-\f{9}{343}p
\mod {p^2}&\t{if $p\e 3\mod 4$,}
\endcases
\\&\sum_{k=0}^{p-1}\f{k^3\b{2k}k^2\b{4k}{2k}}{648^k}
\e \cases \f{1436}{16807}x^2+\f{792}{16807}p\mod {p^2}&\t{if
$p=x^2+4y^2\e 1\mod 4$,}
\\\f {216}{2401}R_1(p)-\f{1510}{16807}p
\mod {p^2}&\t{if $p\e 3\mod 4$.}\endcases
\endalign$$
\endpro
Proof. Taking $a=-\f 14$ and $m=\f {81}8$ in Theorem 3.1 and then
applying (1.1) gives
$$\align &\sum_{k=0}^{p-1}\f{k^2\b{2k}k^2\b{4k}{2k}}{648^k}
\\&\e
\f {16}{49}\sum_{k=0}^{p-1}\f{k\b{2k}k^2\b{4k}{2k}}{648^k}+ \f
6{49}\sum_{k=0}^{p-1}\f{\b{2k}k^2\b{4k}{2k}}{648^k}-\f
3{49}\sum_{k=0}^{p-2}\f{\b{2k}k^2\b{4k}{2k}}{648^k(k+1)}
 \mod {p^3}\endalign$$
 and
$$\align &\sum_{k=0}^{p-1}\f{k^3\b{2k}k^2\b{4k}{2k}}{648^k}
\\&\e
\f {1846}{2401}\sum_{k=0}^{p-1}\f{k\b{2k}k^2\b{4k}{2k}}{648^k}+
\f{435}{2401}\sum_{k=0}^{p-1}\f{\b{2k}k^2\b{4k}{2k}}{648^k}-\f
{144}{2401}\sum_{k=0}^{p-2}\f{\b{2k}k^2\b{4k}{2k}}{648^k(k+1)}
 \mod {p^3}.\endalign$$
 Now applying Theorem 2.10 and (2.21) yields the result.

\section*{4. Congruences for $\sum_{k=0}^{p-2}
\b ak\b{-1-a}k\b{2k}k\f 1{m^k(k+1)^2}$ and $\sum_{k=0}^{p-2} \b
ak\b{-1-a}k\b{2k}k\f 1{m^k(k+1)^3}$
 modulo $p^2$}

 \pro{Theorem 4.1} Let $p$ be an odd prime, $a,m\in\Bbb Z_p$, $m\not\e 0\mod p$
  and
 $a\not \e 0,-1\mod p$. Then
$$\aligned 2a(a+1)\sum_{k=0}^{p-2}\f {\b ak\b{-1-a}k\b{2k}k}{m^k(k+1)^2}
&\e (m-4)\sum_{k=0}^{p-1}\f {k\b ak\b{-1-a}k\b{2k}k}{m^k}
+2\sum_{k=0}^{p-1}\f {\b ak\b{-1-a}k\b{2k}k}{m^k} \\&\q+(4a(a+1)-2)
\sum_{k=0}^{p-2}\f {\b
ak\b{-1-a}k\b{2k}k}{m^k(k+1)}\mod{p^3}\endaligned\tag 4.1$$ and
$$\aligned &2a(a+1)\sum_{k=0}^{p-2}\f {\b ak\b{-1-a}k\b{2k}k}{m^k(k+1)^3}
\e-m+\Big(2m-8-\f{m-4}{a(a+1)}\Big) \sum_{k=0}^{p-1}\f {k\b
ak\b{-1-a}k\b{2k}k}{m^k} \\&\q\qq\qq\qq\qq+\Big(m-\f
2{a(a+1)}\Big)\sum_{k=0}^{p-1}\f {\b ak\b{-1-a}k\b{2k}k}{m^k}
\\&\q\qq\qq\qq\qq+\Big(8a(a+1)-2+\f 2{a(a+1)}\Big) \sum_{k=0}^{p-2}\f {\b
ak\b{-1-a}k\b{2k}k}{m^k(k+1)} \mod {p^3}.
\endaligned\tag 4.2$$
\endpro
Proof. By (3.1),
 $$\align&(k+1)\b a{k+1}\b{-1-a}{k+1}\b{2(k+1)}{k+1}
 \\&=\Big(4k-2+\f{2-4a(a+1)}{k+1}+\f{2a(a+1)}{(k+1)^2}\Big)
 \b ak\b{-1-a}k\b{2k}k\endalign$$
 and so
$$\align&\b a{k+1}\b{-1-a}{k+1}\b{2(k+1)}{k+1}
 \\&=\Big(4-\f 6{k+1}+\f{2-4a(a+1)}{(k+1)^2}+\f{2a(a+1)}{(k+1)^3}\Big)
 \b ak\b{-1-a}k\b{2k}k.\endalign$$
Thus,
$$\align\sum_{k=0}^{p-1}\f{k\b ak\b{-1-a}k\b{2k}k}{m^k}
&=\sum_{k=0}^{p-2}\f{(k+1)\b
a{k+1}\b{-1-a}{k+1}\b{2(k+1)}{k+1}}{m^{k+1}} \\&=\f 1m
\sum_{k=0}^{p-2}\Big(4k-2+\f{2-4a(a+1)}{k+1}+\f{2a(a+1)}{(k+1)^2}\Big)
\f{\b
ak\b{-1-a}k\b{2k}k}{m^k}
\endalign$$
and
$$\align\sum_{k=0}^{p-1}\f{\b ak\b{-1-a}k\b{2k}k}{m^k}
&=1+\sum_{k=0}^{p-2}\f{\b
a{k+1}\b{-1-a}{k+1}\b{2(k+1)}{k+1}}{m^{k+1}}
\\&=1+\f 1m\sum_{k=0}^{p-2}\Big(4-\f 6{k+1}+\f{2-4a(a+1)}{(k+1)^2}+\f{2a(a+1)}{(k+1)^3}\Big)
 \f{\b ak\b{-1-a}k\b{2k}k}{m^k}.
 \endalign$$
 Note that $\b a{p-1}\e \b{-1-a}{p-1}\e \b{2(p-1)}{p-1}\e 0\mod
p$. We then obtain (4.1) and
$$\align 2a(a+1)\sum_{k=0}^{p-2}\f {\b ak\b{-1-a}k\b{2k}k}{m^k(k+1)^3}
&\e -m+(m-4)\sum_{k=0}^{p-1}\f {\b ak\b{-1-a}k\b{2k}k}{m^k}
+6\sum_{k=0}^{p-2}\f {\b ak\b{-1-a}k\b{2k}k}{m^k(k+1)}
\\&\q+(4a(a+1)-2) \sum_{k=0}^{p-2}\f {\b
ak\b{-1-a}k\b{2k}k}{m^k(k+1)^2}\mod{p^3}.\endalign$$ Combining the
two congruences yields (4.2). The proof is now complete.
\par\q
\pro{Corollary 4.1} Let $p$ be an odd prime, $m\in\Bbb Z_p$ and
 $m\not \e 0\mod p$. Then
$$\align&\sum_{k=0}^{\f{p-1}2}\f{\b{2k}k^3}{(16m)^k(k+1)^2}
\e (8-2m)\sum_{k=0}^{\f{p-1}2}\f{k\b{2k}k^3}{(16m)^k}
-4\sum_{k=0}^{\f{p-1}2}\f{\b{2k}k^3}{(16m)^k}
+6\sum_{k=0}^{\f{p-1}2}\f{\b{2k}k^3}{(16m)^k(k+1)}\mod {p^3},
\\&\sum_{k=0}^{(p-1)/2}\f{\b{2k}k^3}{(16m)^k(k+1)^3}
\e 2m-12(m-4)\sum_{k=0}^{(p-1)/2}\f{k\b{2k}k^3}{(16m)^k}
-2(m+8)\sum_{k=0}^{(p-1)/2}\f{\b{2k}k^3}{(16m)^k}
\\&\qq\qq\qq\qq\qq\qq+24\sum_{k=0}^{(p-1)/2}\f{\b{2k}k^3}{(16m)^k(k+1)}\mod {p^3}.
\endalign$$
\endpro
Proof. Taking $a=-\f 12$ in Theorem 4.1 and noting that
$\b{-1/2}k=\b{2k}k(-4)^{-k}$ and $p\mid \b{2k}k$ for $\f p2<k<p$
yields the result. \vskip 0.25cm
 \pro{Theorem 4.2} Let $p$ be an odd prime such
that $p\e 1,2,4\mod 7$ and so $p=x^2+7y^2$. Then
$$\align &\sum_{k=0}^{(p-1)/2}\f{\b{2k}k^3}{(k+1)^2}\e
-68y^2+p\mod {p^2},
\\&\sum_{k=0}^{(p-1)/2}\f{\b{2k}k^3}{(k+1)^3}\e
\f 18-\f{201}2y^2-\f 94p\mod{p^2},
\\&(-1)^{\f{p-1}2}\sum_{k=0}^{(p-1)/2}\f{\b{2k}k^3}{4096^k(k+1)^2}\e
-1136y^2+64p\mod {p^2},
 \\&(-1)^{\f{p-1}2}\sum_{k=0}^{(p-1)/2}\f{\b{2k}k^3}{4096^k(k+1)^3}\e
512(-1)^{\f{p-1}2}+6432y^2-648p\mod{p^2}.
\endalign$$
\endpro
Proof. Taking $m=\f 1{16}$ in Corollary 4.1 and then applying (2.4),
(2.5) and Theorem 2.3 gives
$$\align \sum_{k=0}^{(p-1)/2}\f{\b{2k}k^3}{(k+1)^2}&\e\f{63}8
\sum_{k=0}^{(p-1)/2}k\b{2k}k^3-4\sum_{k=0}^{(p-1)/2}\b{2k}k^3 +6
\sum_{k=0}^{(p-1)/2}\f{\b{2k}k^3}{k+1}
\\&\e \f{63}8\cdot \f 8{21}(28y^2-p)-4(-28y^2+2p)+6(-44y^2+2p)
\\&=-68y^2+p\mod{p^2}
\endalign$$
and
$$\align \sum_{k=0}^{(p-1)/2}\f{\b{2k}k^3}{(k+1)^3}&\e\f 18+\f{189}4
\sum_{k=0}^{(p-1)/2}k\b{2k}k^3-\f{129}8\sum_{k=0}^{(p-1)/2}\b{2k}k^3
+24 \sum_{k=0}^{(p-1)/2}\f{\b{2k}k^3}{k+1}
\\&\e \f 18+\f{189}4\cdot \f 8{21}(28y^2-p)-\f{129}8(-28y^2+2p)+24(-44y^2+2p)
\\&=\f 18-\f{201}2y^2-\f 94p\mod{p^2}.
\endalign$$
Similarly, taking $m=256$ in Corollary 4.1 and then applying (2.6),
(2.7) and Theorem 2.3 yields the remaining congruences. \vskip
0.25cm
\par{\bf Remark 4.1} Let $p$ be an odd prime such that $p\e 1,2,4\mod 7$
and so $p=x^2+7y^2$. In [S10], the author conjectured that
$$\align&\sum_{k=0}^{(p-1)/2}\f{\b{2k}k^3}{(k+1)^2}\e
-68y^2+p-\f{p^2}{4y^2}\mod {p^3},
\\&\sum_{k=0}^{(p-1)/2}\f{\b{2k}k^3}{(k+1)^3}\e
\f 18-\f{201}2y^2-\f 94p-\f{39p^2}{32y^2}\mod {p^3},
\\&(-1)^{\f{p-1}2}\sum_{k=0}^{(p-1)/2}\f{\b{2k}k^3}{4096^k(k+1)^2}\e
-1136y^2+64p+\f{2p^2}{y^2}\mod {p^3},
\\&(-1)^{\f{p-1}2}\sum_{k=0}^{(p-1)/2}\f{\b{2k}k^3}{4096^k(k+1)^3}\e
512(-1)^{\f{p-1}2}+6432y^2-648p-\f{6p^2}{y^2}\mod {p^3}.
\endalign$$
\vskip 0.25cm
 \pro{Theorem 4.3} Suppose that $p$ is a prime of the
form $3k+1$ and so $p=x^2+3y^2$. Then
$$\align&\sum_{k=0}^{(p-1)/2}\f{\b{2k}k^3}{16^k(k+1)^2}\e
-24y^2+2p\mod {p^2},
\\&\sum_{k=0}^{(p-1)/2}\f{\b{2k}k^3}{16^k(k+1)^3}\e
2-24y^2\mod {p^2},
\\&(-1)^{\f{p-1}2}\sum_{k=0}^{(p-1)/2}\f{\b{2k}k^3}{256^k(k+1)^2}\e
-48y^2+8p\mod {p^2},
\\&(-1)^{\f{p-1}2}\sum_{k=0}^{(p-1)/2}\f{\b{2k}k^3}{256^k(k+1)^3}\e
32(-1)^{\f{p-1}2}+96y^2-24p\mod {p^2}.\endalign$$
\endpro
Proof. Taking $m=1,16$ in Corollary 4.1 and then applying
(2.8)-(2.11) and Theorem 2.4 yields the result. \vskip 0.25cm
\pro{Theorem 4.4} Suppose that $p$ is a prime such that $p\e 1,3\mod
8$ and so $p=x^2+2y^2$. Then
$$\align&(-1)^{\f{p-1}2}\sum_{k=0}^{(p-1)/2}\f{\b{2k}k^3}{(-64)^k(k+1)^2}
\e -8y^2\mod {p^2},
\\&(-1)^{\f{p-1}2}\sum_{k=0}^{(p-1)/2}\f{\b{2k}k^3}{(-64)^k(k+1)^3}
\e -8(-1)^{\f{p-1}2}-32y^2+8p\mod {p^2}.\endalign$$
\endpro
Proof. Taking $m=-4$ in Corollary 4.1 and then applying
(2.12)-(2.13) and Theorem 2.5 yields the result.
 \vskip 0.25cm
\pro{Theorem 4.5} Suppose that $p$ is a prime of the form $4k+1$ and
so $p=x^2+4y^2$. Then
$$\align&(-1)^{\f{p-1}4}\sum_{k=0}^{(p-1)/2}\f{\b{2k}k^3}{(-512)^k(k+1)^2}\e
64y^2-8p\mod {p^2},
\\&(-1)^{\f{p-1}4}\sum_{k=0}^{(p-1)/2}\f{\b{2k}k^3}{(-512)^k(k+1)^3}\e
-64(-1)^{\f{p-1}4}-384y^2+72p\mod {p^2}.\endalign$$
\endpro
Proof. Taking $m=-32$ in Corollary 4.1 and then applying
(2.14)-(2.15) and Theorem 2.6 yields the result. \vskip 0.25cm
\pro{Theorem 4.6} Let $p$ be an odd prime. Then

$$\align &\sum_{k=0}^{\f{p-1}2}\f{\b{2k}k^3}{(-8)^k(k+1)^2}\e
\cases -32y^2+p\mod {p^2}&\t{if $p=x^2+4y^2\e 1\mod 4$,}
\\3p+3R_1(p)\mod {p^2}
& \t{if $p\e 3\mod 4$,}
\endcases
\\&\sum_{k=0}^{\f{p-1}2}\f{\b{2k}k^3}{(-8)^k(k+1)^3}\e \cases
-1-48y^2\mod {p^2}&\t{if $p=x^2+4y^2\e 1\mod 4$,}
\\-1+6p+12R_1(p)\mod {p^2}
& \t{if $p\e 3\mod 4$.}
\endcases
\endalign$$
\endpro

Proof. Taking $m=-\f 12$ in Corollary 4.1 gives
$$\align
&\sum_{k=0}^{\f{p-1}2}\f{\b{2k}k^3}{(-8)^k(k+1)^2} \e
9\sum_{k=0}^{\f{p-1}2}\f{k\b{2k}k^3}{(-8)^k}
-4\sum_{k=0}^{\f{p-1}2}\f{\b{2k}k^3}{(-8)^k}+6
\sum_{k=0}^{\f{p-1}2}\f{\b{2k}k^3}{(-8)^k(k+1)}\mod {p^2},
\\&\sum_{k=0}^{\f{p-1}2}\f{\b{2k}k^3}{(-8)^k(k+1)^3} \e
-1+54\sum_{k=0}^{\f{p-1}2}\f{k\b{2k}k^3}{(-8)^k}
-15\sum_{k=0}^{\f{p-1}2}\f{\b{2k}k^3}{(-8)^k}+24
\sum_{k=0}^{\f{p-1}2}\f{\b{2k}k^3}{(-8)^k(k+1)}\mod {p^2}.
\endalign$$
Now applying (3.2), (3.3) and Theorem 2.7 yields the result. \vskip
0.25cm \pro{Theorem 4.7} Let $p>3$ be a prime. Then
$$\align
&\sum_{k=0}^{p-2}\f{\b{2k}k^2\b{3k}k}{(-192)^k(k+1)^2} \e\cases \f
14{x^2}-2p\mod{p^2}\qq\t{if $3\mid p-1$ and so $4p=x^2+27y^2$,}
\\13(2p+1)\b{[2p/3]}{[p/3]}^2+\f 32p\mod {p^2}\qq\qq\q\t{if $3\mid p-2$,}
\endcases
\\&\sum_{k=0}^{p-2}\f{\b{2k}k^2\b{3k}k}{(-192)^k(k+1)^3}
\e\cases -16+\f {51}8x^2-9p\mod {p^2}\q\t{if $3\mid p-1$ and
$4p=x^2+27y^2$,}\\-16+\f{115}2(2p+1)\b{[2p/3]}{[p/3]}^2-\f
{15}{4}p\mod {p^2} \q\  \t{if $3\mid p-2$.}\endcases
\endalign$$
\endpro
Proof. Putting $a=-\f 13$ and $m=-\f{64}9$ in Theorem 4.1 and then
applying (1.1) gives
$$ \sum_{k=0}^{p-2}\f{\b{2k}k^2\b{3k}{k}}{(-192)^k(k+1)^2}\e
25\sum_{k=0}^{p-1}\f{k\b{2k}k^2\b{3k}{k}}{(-192)^k} -\f
92\sum_{k=0}^{p-1}\f{\b{2k}k^2\b{3k}{k}}{(-192)^k} +\f
{13}2\sum_{k=0}^{p-2}\f{\b{2k}k^2\b{3k}{k}}{(-192)^k(k+1)}\mod{p^3}$$
and
$$ \align&\sum_{k=0}^{p-2}\f{\b{2k}k^2\b{3k}{k}}{(-192)^k(k+1)^3}
\\&\e
-16+\f{325}2\sum_{k=0}^{p-1}\f{k\b{2k}k^2\b{3k}{k}}{(-192)^k} -\f
{17}4\sum_{k=0}^{p-1}\f{\b{2k}k^2\b{3k}{k}}{(-192)^k} +\f
{115}4\sum_{k=0}^{p-2}\f{\b{2k}k^2\b{3k}{k}}{(-192)^k(k+1)}\mod{p^3}.
\endalign$$
Now applying Theorem 2.8 and (2.19) yields the result. \vskip 0.25cm
 \pro{Theorem 4.8} Let $p>3$ be a prime. Then
$$\align
&\sum_{k=0}^{p-2}\f{\b{2k}k^2\b{4k}{2k}}{(-144)^k(k+1)^2} \e\cases
-\f{40}3y^2+\f{2p}3\mod {p^2}&\t{if $p=x^2+3y^2\e 1\mod 3$,}
\\\f {88}9R_3(p)+\f {14}9p\mod {p^2}&\t{if $p\e 2\mod 3$,}
\endcases
\\&\sum_{k=0}^{p-2}\f{\b{2k}k^2\b{4k}{2k}}{(-144)^k(k+1)^3}
\e\cases -6-\f{376}{9}y^2+\f{56}{9}p\mod {p^2} &\t{if $p=x^2+3y^2\e
1\mod 3$,}
\\-6+\f {1360}{27}R_3(p)+\f {20}{27}p\mod {p^2}&\t{if $p\e 2\mod 3$.}
\endcases
\endalign$$
\endpro
Proof. Taking $a=-\f 14$ and $m=-\f 94$ in Theorem 4.1 and then
applying (1.1) gives
$$\align &\sum_{k=0}^{p-2}\f{\b{2k}k^2\b{4k}{2k}}{(-144)^k(k+1)^2}
\\&\e \f {50}{3}\sum_{k=0}^{p-1}\f{k\b{2k}k^2\b{4k}{2k}}{(-144)^k}-
\f{16}3\sum_{k=0}^{p-1}\f{\b{2k}k^2\b{4k}{2k}}{(-144)^k}+\f
{22}3\sum_{k=0}^{p-2}\f{\b{2k}k^2\b{4k}{2k}}{(-144)^k(k+1)}
 \mod {p^3}\endalign$$
 and
$$\align &\sum_{k=0}^{p-2}\f{\b{2k}k^2\b{4k}{2k}}{(-144)^k(k+1)^3}
\\&\e -6+\f {1100}{9}\sum_{k=0}^{p-1}\f{k\b{2k}k^2\b{4k}{2k}}{(-144)^k}-
\f{202}9\sum_{k=0}^{p-1}\f{\b{2k}k^2\b{4k}{2k}}{(-144)^k}+\f
{340}9\sum_{k=0}^{p-2}\f{\b{2k}k^2\b{4k}{2k}}{(-144)^k(k+1)}
 \mod {p^3}.\endalign$$
Now applying Theorem 2.9 and (2.20) yields the result. \vskip 0.25cm
\pro{Theorem 4.9} Let $p>3$ be a prime. Then
$$\align
&\sum_{k=0}^{p-2}\f{\b{2k}k^2\b{4k}{2k}}{648^k(k+1)^2} \e\cases \f
{112}9x^2-\f {55}9p\mod {p^2}&\t{if $p=x^2+4y^2\e 1\mod 4$,}
\\-11R_1(p)-\f{1}9 p\mod {p^2}&\t{if $p\e 3\mod 4$,}
\endcases
\\&\sum_{k=0}^{p-2}\f{\b{2k}k^2\b{4k}{2k}}{648^k(k+1)^3} \e\cases
27-\f{740}{27}x^2+\f{248}{27}p\mod {p^2}&\t{if $p=x^2+4y^2\e 1\mod
4$,}
\\27-\f{170}3R_1(p)+\f{122}{27}p\mod {p^2}&\t{if $p\e 3\mod 4$.}
\endcases
\endalign$$
\endpro
Proof. Taking $a=-\f 14$ and $m=\f {81}8$ in Theorem 4.1 and then
applying (1.1) gives
$$\align &\sum_{k=0}^{p-2}\f{\b{2k}k^2\b{4k}{2k}}{648^k(k+1)^2}
\\&\e
-\f {49}{3}\sum_{k=0}^{p-1}\f{k\b{2k}k^2\b{4k}{2k}}{648^k}-
\f{16}3\sum_{k=0}^{p-1}\f{\b{2k}k^2\b{4k}{2k}}{648^k}+\f
{22}3\sum_{k=0}^{p-2}\f{\b{2k}k^2\b{4k}{2k}}{648^k(k+1)}
 \mod {p^3}\endalign$$
 and
$$\align &\sum_{k=0}^{p-2}\f{\b{2k}k^2\b{4k}{2k}}{648^k(k+1)^3}
\\&\e
27-\f {1078}{9}\sum_{k=0}^{p-1}\f{k\b{2k}k^2\b{4k}{2k}}{648^k}-
\f{499}9\sum_{k=0}^{p-1}\f{\b{2k}k^2\b{4k}{2k}}{648^k}+\f
{340}{9}\sum_{k=0}^{p-2}\f{\b{2k}k^2\b{4k}{2k}}{648^k(k+1)}
 \mod {p^3}.\endalign$$
 Now applying Theorem 2.10 and (2.21) yields the result.

\section*{5. Congruences for
$\sum_{k=0}^{p-1}\f{\b ak\b {-1-a}k\b {2k}k}{m^k(2k-1)}$ and
$\sum_{k=0}^{p-1}\f{\b{2k}k^3}{m^k(2k-1)^2}$ modulo $p^2$}

\pro{Theorem 5.1} Suppose that $p$ is an odd prime, $a,m\in\Bbb
Z_p$, $m\not\e 0\mod p$ and $a\not\e 0,-1\mod p$. Then
$$\align\sum_{k=0}^{p-1}\f{\b ak\b{-1-a}k\b{2k}k}{m^k(2k-1)}
&\e \Big(\f 8m-2\Big)\sum_{k=0}^{p-1}\f{k\b ak\b{-1-a}k\b{2k}k}{m^k}
-\sum_{k=0}^{p-1}\f{\b ak\b{-1-a}k\b{2k}k}{m^k}
\\&\q-\f{8a(a+1)}m\sum_{k=0}^{p-2}\f{\b ak\b{-1-a}k\b{2k}k}{m^k(k+1)}\mod {p^3}.
\endalign$$
\endpro
Proof. It is easy to see that
$$\b a{k+1}\b{-1-a}{k+1}\b{2k+2}{k+1}=2(2k+1)\Big(\f 1{k+1}
-\f 1{(k+1)^2}-\f{a(a+1)}{(k+1)^3}\Big)\b ak\b{-1-a}k\b{2k}k.$$
Thus,
$$\align&\sum_{k=1}^{p-1}\f{\b ak\b{-1-a}k\b{2k}k}{m^k(2k-1)}
\\&=\sum_{k=0}^{p-2}\f{\b
a{k+1}\b{-1-a}{k+1}\b{2k+2}{k+1}}{m^{k+1}(2(k+1)-1)} =\f
2m\sum_{k=0}^{p-2}\Big(\f 1{k+1}-\f 1{(k+1)^2}-\f
{a(a+1)}{(k+1)^3}\Big)\f{\b ak\b{-1-a}k\b{2k}k}{m^k}.\endalign$$
Now, applying Theorem 4.1 we get
$$\align&\sum_{k=1}^{p-1}\f{\b ak\b{-1-a}k\b{2k}k}{m^k(2k-1)}
\\&\e\f 2m\sum_{k=0}^{p-2}\f{\b ak\b{-1-a}k\b{2k}k}{m^k(k+1)}
-\f 1{a(a+1)m}\Big((m-4)\sum_{k=0}^{p-1}\f{k\b
ak\b{-1-a}k\b{2k}k}{m^k} +2\sum_{k=0}^{p-1}\f {\b
ak\b{-1-a}k\b{2k}k}{m^k} \\&\q+(4a(a+1)-2) \sum_{k=0}^{p-2}\f {\b
ak\b{-1-a}k\b{2k}k}{m^k(k+1)}\Big) -\f
1m\Big(-m+\Big(2m-8-\f{m-4}{a(a+1)}\Big)\\&\q\times
\sum_{k=0}^{p-1}\f {k\b ak\b{-1-a}k\b{2k}k}{m^k}+\Big(m-\f
2{a(a+1)}\Big)\sum_{k=0}^{p-1}\f {\b ak\b{-1-a}k\b{2k}k}{m^k}
\\&\q+\Big(8a(a+1)-2+\f 2{a(a+1)}\Big) \sum_{k=0}^{p-2}\f {\b
ak\b{-1-a}k\b{2k}k}{m^k(k+1)}\Big)\mod{p^3},\endalign$$ which yields
the result.
\par\q
\pro{Corollary 5.1}  Let $p$ be an odd prime, $m\in\Bbb Z_p$ and
$m\not\e 0\mod p$. Then
$$\align\sum_{k=0}^{p-1}\f{\b{2k}k^3}{(16m)^k(2k-1)}
&\e \Big(\f 8m-2\Big)\sum_{k=0}^{(p-1)/2}
\f{k\b{2k}k^3}{(16m)^k}-\sum_{k=0}^{(p-1)/2} \f{\b{2k}k^3}{(16m)^k}
\\&\q+\f{2}m\sum_{k=0}^{(p-1)/2} \f{\b{2k}k^3}{(16m)^k(k+1)} \mod
{p^3}.\endalign$$\endpro
 Proof. Taking $a=-\f 12$ in Theorem 5.1
yields the result.
\par\q
\pro{Theorem 5.2}  Let $p$ be an odd prime, $m\in\Bbb Z_p$ and
$m\not\e 0\mod p$. Then
$$\align\sum_{k=0}^{p-1}\f{\b{2k}k^3}{(16m)^k(2k-1)^2}
&\e \Big(4-\f {16}m\Big)\sum_{k=0}^{(p-1)/2}
\f{k\b{2k}k^3}{(16m)^k}+\Big(1+\f 4m\Big)\sum_{k=0}^{(p-1)/2}
\f{\b{2k}k^3}{(16m)^k}\\&\q-\f 6m\sum_{k=0}^{(p-1)/2}
\f{\b{2k}k^3}{(16m)^k(k+1)}\mod {p^3}.\endalign$$
\endpro
Proof. Since
$$\f{\b{2(k+1)}{k+1}^3}{(2k+1)^2}
=\f{(2\f{2k+1}{k+1}\b{2k}k)^3}{(2k+1)^2}=8\Big(\f 2{(k+1)^2}-\f
1{(k+1)^3}\Big)\b{2k}k^3,$$ applying Corollary 4.1 we see that
$$\align&\sum_{k=1}^{p-1}\f{\b{2k}k^3}{(16m)^k(2k-1)^2}
\\&=\sum_{k=0}^{p-2}\f{\b{2k+2}{k+1}^3}{(16m)^{k+1}(2(k+1)-1)^2}
=\f 1{2m}\sum_{k=0}^{p-2}\Big(\f 2{(k+1)^2}-\f
1{(k+1)^3}\Big)\f{\b{2k}k^3}{(16m)^k}
\\&=\f 1{2m}\Big(2(8-2m)\sum_{k=0}^{(p-1)/2}
\f{k\b{2k}k^3}{(16m)^k}-8\sum_{k=0}^{(p-1)/2} \f{\b{2k}k^3}{(16m)^k}
+12\sum_{k=0}^{(p-1)/2} \f{\b{2k}k^3}{(16m)^k(k+1)}
\\&\q-2m+12(m-4)\sum_{k=0}^{(p-1)/2}
\f{k\b{2k}k^3}{(16m)^k}+2(m+8)\sum_{k=0}^{(p-1)/2}
\f{\b{2k}k^3}{(16m)^k}\\&\q-24\sum_{k=0}^{(p-1)/2}
\f{\b{2k}k^3}{(16m)^k(k+1)}\Big)\mod {p^3},
\endalign$$
which yields the result.
\par\q
 \pro{Theorem 5.3} Let $p$ be an odd prime
such that $p\e 1,2,4\mod 7$ and so $p=x^2+7y^2$. Then
$$\align&\sum_{k=0}^{p-1}\f{\b{2k}k^3}{2k-1}\e
-36y^2+14p\mod {p^2},
\\&\sum_{k=0}^{p-1}\f{\b{2k}k^3}{(2k-1)^2}\e  -284y^2+34p\mod {p^2},
\\&(-1)^{\f{p-1}2}\sum_{k=0}^{p-1}\f{\b{2k}k^3}{4096^k(2k-1)}
\e 22y^2-\f 74p\mod {p^2},
\\&(-1)^{\f{p-1}2}\sum_{k=0}^{p-1}\f{\b{2k}k^3}{4096^k(2k-1)^2}
\e -17y^2+\f{97}{64}p\mod {p^2}.\endalign$$
\endpro
Proof. Taking $m=\f 1{16}$ in Corollary 5.1 and then applying (2.4),
(2.5) and Theorem 2.3 gives
$$\align\sum_{k=0}^{p-1}\f{\b{2k}k^3}{2k-1}
&\e 126\sum_{k=0}^{(p-1)/2}k\b{2k}k^3-\sum_{k=0}^{(p-1)/2}\b{2k}k^3
+32\sum_{k=0}^{(p-1)/2}\f{\b{2k}k^3}{k+1}
\\&\e 126\cdot\f 8{21}(28y^2-p)-(2p-28y^2)+32(-44y^2+2p)
\\&=-36y^2+14p\mod {p^2},\endalign$$
and putting $m=\f 1{16}$ in Theorem 5.2 and then applying (2.4),
(2.5) and Theorem 2.3 gives
$$\align\sum_{k=0}^{p-1}\f{\b{2k}k^3}{(2k-1)^2}
&\e
-252\sum_{k=0}^{(p-1)/2}k\b{2k}k^3+65\sum_{k=0}^{(p-1)/2}\b{2k}k^3
-96\sum_{k=0}^{(p-1)/2}\f{\b{2k}k^3}{k+1}
\\&\e -252\cdot\f 8{21}(28y^2-p)+65(2p-28y^2)-96(-44y^2+2p)
\\&=-284y^2+34p\mod {p^2}.\endalign$$
On the other hand, taking $m=256$ in Corollary 5.1 and then applying
(2.6), (2.7) and Theorem 2.3 gives
$$\align&\sum_{k=0}^{p-1}\f{\b{2k}k^3}{4096^k(2k-1)}
\\&\e -\f{63}{32}\sum_{k=0}^{(p-1)/2}\f{k\b{2k}k^3}{4096^k}
-\sum_{k=0}^{(p-1)/2}\f{\b{2k}k^3}{4096^k}
+\f{1}{128}\sum_{k=0}^{(p-1)/2}\f{\b{2k}k^3}{4096^k(k+1)}
\\&\e -\f{63}{32}\cdot \f 5{42}(-1)^{\f{p-1}2}(28y^2-p)-(-1)^{\f{p-1}2}
(-28y^2+2p)+\f 1{128}(-1)^{\f{p-1}2}(72y^2+2p)
\\&=(-1)^{\f{p-1}2}\big(22y^2-\f 74p\big)\mod
{p^2}.\endalign$$ Putting $m=256$ in Theorem 5.2 and then applying
(2.6), (2.7) and Theorem 2.3 gives
$$\align&\sum_{k=0}^{p-1}\f{\b{2k}k^3}{4096^k(2k-1)^2}
\\&\e -12\sum_{k=0}^{(p-1)/2}\f{k\b{2k}k^3}{4096^k}
+\f{65}{64}\sum_{k=0}^{(p-1)/2}\f{\b{2k}k^3}{4096^k}
-\f{3}{128}\sum_{k=0}^{(p-1)/2}\f{\b{2k}k^3}{4096^k(k+1)}
\\&\e -12\cdot \f 5{42}(-1)^{\f{p-1}2}(28y^2-p)+\f{65}{64}(-1)^{\f{p-1}2}
(-28y^2+2p)-\f 3{128}(-1)^{\f{p-1}2}(72y^2+2p)
\\&=(-1)^{\f{p-1}2}\big(-17y^2+\f{97}{64}p\big)\mod
{p^2}.\endalign$$ This completes the proof.
\par\q
\pro{Theorem 5.4} Let $p$ be a prime of the form $3k+1$ and so
$p=x^2+3y^2$. Then
$$\align&\sum_{k=0}^{p-1}\f{\b{2k}k^3}{16^k(2k-1)}
\e 4y^2\mod {p^2},
\\&\sum_{k=0}^{p-1}\f{\b{2k}k^3}{16^k(2k-1)^2}\e
-12y^2+2p\mod{p^2},
\\&(-1)^{\f{p-1}2}\sum_{k=0}^{p-1}\f{\b{2k}k^3}{256^k(2k-1)}
\e 8y^2-\f 32p\mod {p^2},
\\&(-1)^{\f{p-1}2}\sum_{k=0}^{p-1}\f{\b{2k}k^3}{256^k(2k-1)^2}
\e -6y^2+\f 54p\mod {p^2}.
\endalign$$
\endpro
Proof. Putting $m=1$ in Corollary 5.1 and Theorem 5.2 and then
applying (2.8), (2.9) and Theorem 2.4 gives the first two
congruences, and taking $m=16$ in Corollary 5.1 and Theorem 5.2 and
then applying (2.10), (2.11) and Theorem 2.4 yields the remaining
congruences.
\par\q
\pro{Theorem 5.5} Let $p$ be a prime such that $p\e 1,3\mod 8$ and
so $p=x^2+2y^2$. Then
$$\align&(-1)^{\f{p-1}2}\sum_{k=0}^{p-1}\f{\b{2k}k^3}{(-64)^k(2k-1)}
\e p-3x^2\mod {p^2},
\\&(-1)^{\f{p-1}2}\sum_{k=0}^{p-1}\f{\b{2k}k^3}{(-64)^k(2k-1)^2}
\e x^2\mod {p^2}.\endalign$$
\endpro
Proof. Putting $m=-4$ in Corollary 5.1 and Theorem 5.2 and then
applying (2.12), (2.13) and Theorem 2.5 yields the result.
\par\q
\pro{Theorem 5.6} Let $p$ be a prime of the form $4k+1$ and so
$p=x^2+4y^2$. Then
$$\align&(-1)^{\f{p-1}4}\sum_{k=0}^{p-1}\f{\b{2k}k^3}{(-512)^k(2k-1)}
\e -3x^2+\f 54p\mod {p^2},
\\&(-1)^{\f{p-1}4}\sum_{k=0}^{p-1}\f{\b{2k}k^3}{(-512)^k(2k-1)^2}
\e  2x^2-\f 58p\mod {p^2}.
\endalign$$
\endpro
Proof. Taking $m=-32$ in Corollary 5.1 and Theorem 5.2 and then
applying (2.14), (2.15) and Theorem 2.6 yields the result. \vskip
0.25cm
 \pro{Theorem 5.7} Let $p$ be an odd prime. Then
$$\align&\sum_{k=0}^{p-1}\f{\b{2k}k^3}{(-8)^k(2k-1)}
\e \cases -4x^2\mod {p^2}&\t{if $p=x^2+4y^2\e 1\mod 4$,}
\\2p-2R_1(p)\mod {p^2}
& \t{if $p\e 3\mod 4$,}\endcases
\\&\sum_{k=0}^{p-1}\f{\b{2k}k^3}{(-8)^k(2k-1)^2} \e
\cases -4x^2+2p\mod {p^2}&\t{if $p=x^2+4y^2\e 1\mod 4$,}
\\6R_1(p)
\mod {p^2} &\t{if $p\e 3\mod 4$.}\endcases
\endalign$$
\endpro
Proof. Taking $m=-\f 12$ in Corollary 5.1 and Theorem 5.2 gives
$$\align&\sum_{k=0}^{p-1}\f{\b{2k}k^3}{(-8)^k(2k-1)}
\e -18\sum_{k=0}^{\f{p-1}2}\f{k\b{2k}k^3}{(-8)^k}-
\sum_{k=0}^{\f{p-1}2}\f{\b{2k}k^3}{(-8)^k}-4
\sum_{k=0}^{\f{p-1}2}\f{\b{2k}k^3}{(-8)^k(k+1)}\mod {p^2},
\\&\sum_{k=0}^{p-1}\f{\b{2k}k^3}{(-8)^k(2k-1)^2}
\e 36\sum_{k=0}^{\f{p-1}2}\f{k\b{2k}k^3}{(-8)^k}-
7\sum_{k=0}^{\f{p-1}2}\f{\b{2k}k^3}{(-8)^k}+12
\sum_{k=0}^{\f{p-1}2}\f{\b{2k}k^3}{(-8)^k(k+1)}\mod
{p^2}.\endalign$$
 Now applying (3.2), (3.3) and Theorem 2.7 yields
the result. \vskip 0.25cm

\pro{Theorem 5.8} Let $p>3$ be a prime. Then
$$ \sum_{k=0}^{p-1}\f{\b{2k}k^2\b{3k}{k}}{(-192)^k(2k-1)}\e\cases -\f
34x^2+\f {9p}8\mod{p^2}\qq\t{if $3\mid p-1$ and so $4p=x^2+27y^2$,}
\\-\f 12(2p+1)\b{[2p/3]}{[p/3]}^2+\f 38p\mod {p^2}\q\ \q\t{if $p\e 2\mod{3}$.}
\endcases$$
\endpro

Proof. Putting $a=-\f 13$ and $m=-\f{64}9$ in Theorem 5.1 and then
applying (1.1) gives
$$ \sum_{k=0}^{p-1}\f{\b{2k}k^2\b{3k}{k}}{(-192)^k(2k-1)}\e
-\f{25}8\sum_{k=0}^{p-1}\f{k\b{2k}k^2\b{3k}{k}}{(-192)^k}
-\sum_{k=0}^{p-1}\f{\b{2k}k^2\b{3k}{k}}{(-192)^k} -\f
14\sum_{k=0}^{p-2}\f{\b{2k}k^2\b{3k}{k}}{(-192)^k(k+1)}\mod{p^3}.$$
Now applying Theorem 2.8 and (2.19) yields the result. \vskip 0.25cm
\pro{Theorem 5.9} Let $p>3$ be a prime. Then
$$\sum_{k=0}^{p-1}\f{\b{2k}k^2\b{4k}{2k}}{(-144)^k(2k-1)} \e\cases
-\f {28}9x^2+\f 89p\mod {p^2}\qq\t{if $p=x^2+3y^2\e 1\mod 3$,}
\\-\f 89R_3(p)+\f 23p\mod {p^2}\q\t{if $p\e 2\mod 3$.}
\endcases$$
\endpro
Proof. Taking $a=-\f 14$ and $m=-\f 94$ in Theorem 5.1 and then
applying (1.1) gives
$$\align &\sum_{k=0}^{p-1}\f{\b{2k}k^2\b{4k}{2k}}{(-144)^k(2k-1)}
\\&\e -\f {50}{9}\sum_{k=0}^{p-1}\f{k\b{2k}k^2\b{4k}{2k}}{(-144)^k}-
\sum_{k=0}^{p-1}\f{\b{2k}k^2\b{4k}{2k}}{(-144)^k}-\f
23\sum_{k=0}^{p-2}\f{\b{2k}k^2\b{4k}{2k}}{(-144)^k(k+1)}
 \mod {p^3}.\endalign$$
 Now applying Theorem 2.9 and (2.20) yields the result.
\vskip 0.25cm

 \pro{Theorem 5.10} Let $p>3$ be a prime. Then
$$\sum_{k=0}^{p-1}\f{\b{2k}k^2\b{4k}{2k}}{648^k(2k-1)}
\e \cases -\f{76}{27}x^2+\f{104}{81}p\mod {p^2}&\t{if $p=x^2+4y^2\e
1\mod 4$,}
\\-\f 29R_1(p)+\f{10}{81}p
\mod {p^2}&\t{if $p\e 3\mod 4$.}
\endcases$$
\endpro
Proof. Taking $a=-\f 14$ and $m=\f {81}8$ in Theorem 5.1 and then
applying (1.1) gives
$$\align &\sum_{k=0}^{p-1}\f{\b{2k}k^2\b{4k}{2k}}{648^k(2k-1)}
\\&\e
-\f {98}{81}\sum_{k=0}^{p-1}\f{k\b{2k}k^2\b{4k}{2k}}{648^k}-
\sum_{k=0}^{p-1}\f{\b{2k}k^2\b{4k}{2k}}{648^k}+\f
4{27}\sum_{k=0}^{p-2}\f{\b{2k}k^2\b{4k}{2k}}{648^k(k+1)}
 \mod {p^3}.\endalign$$
 Now applying Theorem 2.10 and (2.21) yields the result.

\section*{6. Congruences for $\sum_{k=0}^{p-3}\b ak\b{-1-a}k\b{2k}k\f 1
{m^k(k+2)}$ modulo $p^2$}
 \pro{Theorem 6.1} Suppose that $p>3$ is a prime,
$a,m\in\Bbb Z_p$, $m\not\e 0\mod p$ and $a\not\e 0,\pm 1,-2\mod p$.
Then
$$\align&\sum_{k=0}^{p-3}\f{\b ak\b{-1-a}k\b{2k}k}{m^k(k+2)}
\\&\e \f{4-m}{6(a-1)(a+2)}\sum_{k=0}^{p-1}\f{k\b
ak\b{-1-a}k\b{2k}k}{m^k} +\f{m-6}{6(a-1)(a+2)}\sum_{k=0}^{p-1}\f{\b
ak\b{-1-a}k\b{2k}k}{m^k}
\\&\q+\f{2a(a+1)-m}{6(a-1)(a+2)}\sum_{k=0}^{p-2}\f{\b ak\b{-1-a}k\b{2k}k}{m^k(k+1)}\mod {p^3}.
\endalign$$
\endpro
Proof. By (3.1),
$$\align &\f{(k+1)^3\b a{k+1}\b{-1-a}{k+1}\b{2k+2}{k+1}}{k+2}
\\&=\f{(4k^2+2k-4a(a+1))(k+1)+2a(a+1)}{k+2}\b
ak\b{-1-a}k\b{2k}k
\\&=\Big(4k^2-2k+6-4a(a+1)+\f{6(a-1)(a+2)}{k+2}\Big)\b ak\b{-1-a}k\b{2k}k.
\endalign$$
Thus,
$$\align &\sum_{k=0}^{p-2}\f{k^3\b ak\b{-1-a}k\b{2k}k}{m^k(k+1)}
\\&=\sum_{k=0}^{p-3}\f{(k+1)^3\b
a{k+1}\b{-1-a}{k+1}\b{2k+2}{k+1}}{m^{k+1}(k+2)}
\\&=\f 1m\sum_{k=0}^{p-3}\Big(4k^2-2k+6-4a(a+1)+\f{6(a-1)(a+2)}{k+2}\Big)
\f{\b ak\b{-1-a}k\b{2k}k}{m^k}
\endalign$$
On the other hand, since $\f{k^3}{k+1}=k^2-k+1-\f 1{k+1}$ we see
that
$$\sum_{k=0}^{p-2}\f{k^3\b ak\b{-1-a}k\b{2k}k}{m^k(k+1)}
=\sum_{k=0}^{p-2}\Big(k^2-k+1-\f 1{k+1}\Big) \f{\b
ak\b{-1-a}k\b{2k}k}{m^k}.$$ Note that $\b ak\e \b{-1-a}k\e \b{2k}k\e
0\mod p$ for $k=p-1,p-2$. We then get
$$\align &6(a-1)(a+2)\sum_{k=0}^{p-3}
\f{\b ak\b{-1-a}k\b{2k}k}{m^k(k+2)} \\&\e
(m-4)\sum_{k=0}^{p-1}\f{k^2\b ak\b{-1-a}k\b{2k}k}{m^k}
+(2-m)\sum_{k=0}^{p-1}\f{k\b ak\b{-1-a}k\b{2k}k}{m^k}
\\&\q+(4a(a+1)-6+m)\sum_{k=0}^{p-1}
\f{\b ak\b{-1-a}k\b{2k}k}{m^k}-m\sum_{k=0}^{p-2}\f{\b
ak\b{-1-a}k\b{2k}k}{m^k(k+1)}\mod{p^3}.\endalign$$ Now applying the
first congruence in Theorem 3.1 yields the result. \vskip 0.25cm
\pro{Corollary 6.1} Let $p>3$ be a prime, $m\in\Bbb Z_p$ and
$m\not\e 0\mod p$. Then
$$\align&\sum_{k=0}^{(p-1)/2}\f{\b{2k}k^3}{(16m)^k(k+2)}
\\&\e \f{2m-8}{27}\sum_{k=0}^{\f{p-1}2}\f{k\b{2k}k^3}{(16m)^k}
-\f{2m-12}{27}\sum_{k=0}^{\f{p-1}2}\f{\b{2k}k^3}{(16m)^k}
+\f{2m+1}{27} \sum_{k=0}^{\f{p-1}2}\f{\b{2k}k^3}{(16m)^k(k+1)}\mod
{p^3}.\endalign$$ \endpro
Proof. Taking $a=-\f 12$ in Theorem 6.1
gives the result.
\vskip 0.25cm \pro{Theorem 6.2}  Let $p>3$ be a
prime such that $p\e 1,2,4\mod 7$ and so $p=x^2+7y^2$. Then
$$\align &\sum_{k=0}^{(p-1)/2}\f{\b{2k}k^3}{k+2}
 \e -\f{466}{27}y^2+\f{29}{27}p\mod{p^2},
 \\&(-1)^{\f{p-1}2}\sum_{k=0}^{(p-1)/2}\f{\b{2k}k^3}{4096^k(k+2)}
 \e \f{52616}{27}y^2-\f{34}{27}p\mod {p^2}.
 \endalign$$
 \endpro
 Proof. Taking $m=\f 1{16}$ in Corollary 6.1 and then applying
 (2.4), (2.5) and Theorem 2.3 gives
 $$\align \sum_{k=0}^{(p-1)/2}\f{\b{2k}k^3}{k+2}
 &\e -\f 7{24}\sum_{k=0}^{(p-1)/2}k\b{2k}k^3
 +\f{95}{216}\sum_{k=0}^{(p-1)/2}\b{2k}k^3+\f 1{24}
 \sum_{k=0}^{(p-1)/2}\f{\b{2k}k^3}{k+1}
 \\&\e -\f 7{24}\cdot\f 8{21}(28y^2-p)+\f{95}{216}(-28y^2+2p)+
 \f 1{24}(-44y^2+2p)
 \\&=-\f{466}{27}y^2+\f{29}{27}p\mod{p^2},\endalign$$
and taking $m=256$ in Corollary 6.1 and then applying
 (2.6), (2.7) and Theorem 2.3 gives
 $$\align &\sum_{k=0}^{(p-1)/2}\f{\b{2k}k^3}{4096^k(k+2)}
 \\&\e \f {56}3\sum_{k=0}^{(p-1)/2}\f{k\b{2k}k^3}{4096^k}
 -\f{500}{27}\sum_{k=0}^{(p-1)/2}\f{\b{2k}k^3}{4096^k}+19
 \sum_{k=0}^{(p-1)/2}\f{\b{2k}k^3}{4096^k(k+1)}
 \\&\e \f {56}{3}\cdot\f 5{42}(-1)^{\f{p-1}2}(28y^2-p)-\f{500}{27}
 (-1)^{\f{p-1}2}(-28y^2+2p)+
 19(-1)^{\f{p-1}2}(72y^2+2p)
 \\&=(-1)^{\f{p-1}2}\big(\f{52616}{27}y^2-\f{34}{27}p\big)\mod {p^2}.\endalign$$
This completes the proof.
\vskip 0.25cm
\par{\bf Remark 6.1} Let $p$ be an odd prime such that $p\e 1,2,4\mod 7$
and so $p=x^2+7y^2$. In [S10], the author conjectured that
$$\align &\sum_{k=0}^{(p-1)/2}\f{\b{2k}k^3}{k+2}
 \e -\f{466}{27}y^2+\f{29}{27}p+\f{17p^2}{864y^2}\mod {p^3},
\\&\sum_{k=0}^{(p-1)/2}\f{\b{2k}k^3}{k+3}
 \e -\f{36052}{3375}y^2+\f{2378}{3375}p+\f{1421p^2}{108000y^2}\mod
 {p^3},
\\&(-1)^{\f{p-1}2}\sum_{k=0}^{(p-1)/2}\f{\b{2k}k^3}{4096^k(k+2)}
 \e \f{52616}{27}y^2-\f{34}{27}p-\f{20p^2}{27y^2}\mod {p^3},
 \\&(-1)^{\f{p-1}2}\sum_{k=0}^{(p-1)/2}\f{\b{2k}k^3}{4096^k(k+3)}
 \e \f{217125848}{3375}y^2-\f{250882}{3375}p
 -\f{83972p^2}{3375y^2}\mod {p^3}.\endalign$$
\vskip 0.25cm \pro{Theorem 6.3}  Let $p$ be a prime of the form
$3k+1$ and so $p=x^2+3y^2$. Then
$$\align &\sum_{k=0}^{(p-1)/2}\f{\b{2k}k^3}{16^k(k+2)}
 \e \f{64}{27}x^2-\f 43p\mod{p^2},
 \\&(-1)^{\f{p-1}2}\sum_{k=0}^{(p-1)/2}\f{\b{2k}k^3}{256^k(k+2)}
 \e -\f{8}{27}x^2+\f {10}9p\mod {p^2}.
 \endalign$$
 \endpro
Proof. Putting $m=1$ in Corollary 6.1 and then applying
 (2.8), (2.9) and Theorem 2.4 gives the first congruence,
 and taking $m=16$ in Corollary 6.1 and then applying
 (2.10), (2.11) and Theorem 2.4 yields the second congruence.
\vskip 0.25cm
 \pro{Theorem 6.4}  Let $p>3$ be a prime such that
 $p\e 1,3\mod 8$ and so $p=x^2+2y^2$. Then
$$(-1)^{\f{p-1}2}\sum_{k=0}^{(p-1)/2}\f{\b{2k}k^3}{(-64)^k(k+2)}
 \e  2x^2-\f 89p\mod {p^2}.$$
 \endpro
Proof. Taking $m=-4$ in Corollary 6.1 and then applying (2.12),
(2.13) and Theorem 2.5 yields
$$\align &
\sum_{k=0}^{(p-1)/2}\f{\b{2k}k^3}{(-64)^k(k+2)}
\\&\e -\f{16}{27}\sum_{k=0}^{(p-1)/2}\f{k\b{2k}k^3}{(-64)^k}
+\f{20}{27}\sum_{k=0}^{(p-1)/2}\f{\b{2k}k^3}{(-64)^k} -\f 7{27}
\sum_{k=0}^{(p-1)/2}\f{\b{2k}k^3}{(-64)^k(k+1)}
\\&\e -\f{16}{27}\cdot \f 14(-1)^{\f{p-1}2}(8y^2-p)
+\f{20}{27}(-1)^{\f{p-1}2}(-8y^2+2p)-\f
7{27}(-1)^{\f{p-1}2}(-12y^2+2p)
\\&=(-1)^{\f{p-1}2}\big(2x^2-\f 89p\big)\mod {p^2}.
\endalign$$
This proves the theorem. \vskip 0.25cm

\vskip 0.25cm \pro{Theorem 6.5}  Let $p$ be a prime of the form
$4k+1$ and so $p=x^2+4y^2$. Then
$$(-1)^{\f{p-1}4}\sum_{k=0}^{(p-1)/2}\f{\b{2k}k^3}{(-512)^k(k+2)}
 \e -\f{152}{27}x^2+\f{190}{27}p\mod {p^2}.$$
 \endpro
Proof. Taking $m=-32$ in Corollary 6.1 and then applying
 (2.14), (2.15) and Theorem 2.6 yields the result.
 \vskip 0.25cm
\pro{Theorem 6.6} Let $p>3$ be a prime. Then
$$\sum_{k=0}^{(p-1)/2}\f{\b{2k}k^3}{(-8)^k(k+2)}
 \e \cases\f{64}{27}x^2-\f{35}{27}p\mod{p^2}&\t{if $p=x^2+4y^2\e 1\mod
 4$,}\\\f p9\mod {p^2}&\t{if $p\e 3\mod 4$.}
 \endcases$$
 \endpro
Proof. Putting $m=-\f 12$ in Corollary 6.1 and then applying (3.2)
and (3.3) gives
$$\align\sum_{k=0}^{(p-1)/2}\f{\b{2k}k^3}{(-8)^k(k+2)}
 &\e-\f 13\sum_{k=0}^{(p-1)/2}\f{k\b{2k}k^3}{(-8)^k}
 +\f{13}{27}\sum_{k=0}^{(p-1)/2}\f{\b{2k}k^3}{(-8)^k}
\\&\e\cases -\f 13(p-\f 43x^2)+\f{13}{27}(4x^2-2p)=
\f{64}{27}x^2-\f{35}{27}p\mod{p^2}\\\qq\qq\qq\qq\qq\;\t{if
$p=x^2+4y^2\e 1\mod
 4$,}\\-\f 13(-\f 13p)=\f p9\mod {p^2}\q\qq\qq\t{if $p\e 3\mod 4$.}
 \endcases\endalign$$
 This proves the theorem.
 \vskip 0.25cm
 \par{\bf Remark 6.2}
  Let $r\in\{2,3,4,\ldots\}$, and let $p$ be a prime
 such that $p\ge r+1$. Using the method in the proofs of
 Theorems 4.1 and 6.1,
 one may deduce similar congruences modulo $p^2$ for
 $$\sum_{k=0}^{p-1-r}\f{\b{2k}k^3}{m^k(k+r)},\q
\sum_{k=0}^{p-1-r}\f{\b{2k}k^3}{m^k(k+r)^2}\qtq{and}
\sum_{k=0}^{p-1-r}\f{\b{2k}k^3}{m^k(k+r)^3},$$ where $m\in\Bbb Z_p$
and $m\not\e 0\mod p$.
 \pro{Theorem 6.7} Let $p>5$ be a prime. Then
 $$\sum_{k=0}^{p-3}\f{\b{2k}k^2\b{3k}k}{(-192)^k(k+2)} \e \cases
 \f 25x^2-\f
7{15}p\mod {p^3}\qq\t{if $3\mid p-1$ and $4p=x^2+27y^2$,}
\\-(2p+1)\b{[2p/3]}{[p/3]}^2-\f p3\mod {p^2}\q\qq\q\ \t{if $3\mid p-2$.}
\endcases$$
\endpro
Proof. Taking $a=-\f 13$ and $m=-\f{64}9$ in Theorem 6.1 and then
applying (1.1) gives
$$\align &\sum_{k=0}^{p-3}\f{\b{2k}k^2\b{3k}k}{(-192)^k(k+2)}
\\&\e
-\f 56\sum_{k=0}^{p-1}\f{k\b{2k}k^2\b{3k}k}{(-192)^k}+\f{59}{60}
\sum_{k=0}^{p-1}\f{\b{2k}k^2\b{3k}k}{(-192)^k}-\f
12\sum_{k=0}^{p-2}\f{\b{2k}k^2\b{3k}k}{(-192)^k(k+1)}
 \mod {p^3}.\endalign$$
 Now applying Theorem 2.8 and (2.19) yields the result.

 \pro{Theorem 6.8} Let $p>7$ be a prime. Then
$$\sum_{k=0}^{p-3}\f{\b{2k}k^2\b{4k}{2k}}{(-144)^k(k+2)}
\e\cases -\f{32}5y^2+\f{16}{15}p\mod {p^2} &\t{if $p=x^2+3y^2\e
1\mod 3$,}
\\-\f 4{21}R_3(p)\mod {p^2}&\t{if $p\e 2\mod 3$.}
\endcases$$
\endpro
Proof. Taking $a=-\f 14$ and $m=-\f 94$ in Theorem 6.1 and then
applying (1.1) gives
$$\align &\sum_{k=0}^{p-3}\f{\b{2k}k^2\b{4k}{2k}}{(-144)^k(k+2)}
\\&\e
-\f
{10}{21}\sum_{k=0}^{p-1}\f{k\b{2k}k^2\b{4k}{2k}}{(-144)^k}+\f{22}{35}
\sum_{k=0}^{p-1}\f{\b{2k}k^2\b{4k}{2k}}{(-144)^k}-\f
17\sum_{k=0}^{p-2}\f{\b{2k}k^2\b{4k}{2k}}{(-144)^k(k+1)}
 \mod {p^3}.\endalign$$
 Now applying Theorem 2.9 and (2.20) yields the result.

\pro{Theorem 6.9} Let $p>7$ be a prime. Then
$$\sum_{k=0}^{p-3}\f{\b{2k}k^2\b{4k}{2k}}{648^k(k+2)}
\e\cases \f 87x^2-\f 5{21}p\mod {p^2} &\t{if $p=x^2+4y^2\e 1\mod
4$,}
\\-\f p3-\f 65R_1(p)\mod {p^2}&\t{if $p\e 3\mod 4$.}
\endcases$$
\endpro
Proof. Taking $a=-\f 14$ and $m=\f {81}8$ in Theorem 6.1 and then
applying (1.1) gives
$$\align &\sum_{k=0}^{p-3}\f{\b{2k}k^2\b{4k}{2k}}{648^k(k+2)}
\\&\e
\f {7}{15}\sum_{k=0}^{p-1}\f{k\b{2k}k^2\b{4k}{2k}}{648^k}-\f{11}{35}
\sum_{k=0}^{p-1}\f{\b{2k}k^2\b{4k}{2k}}{648^k}+\f
45\sum_{k=0}^{p-2}\f{\b{2k}k^2\b{4k}{2k}}{648^k(k+1)}
 \mod {p^3}.\endalign$$
 Now applying Theorem 2.10 and (2.21) yields the result.

\end{document}